\documentclass[a4paper,11pt]{article}
\usepackage{amsmath,amssymb,graphicx}
\usepackage{subfigure}
\usepackage{lscape}
\usepackage{mathrsfs}
\usepackage{multirow}
\usepackage{verbatim}
\usepackage{color}
\usepackage{ntheorem}
\usepackage{mathtools}
\usepackage{multirow}
\usepackage{xcolor}
\usepackage{graphicx}
\usepackage{epstopdf}
\usepackage{enumerate}
\usepackage{amsfonts}
\usepackage{epsfig}
\usepackage{amscd}
\usepackage{stmaryrd}
\usepackage{arydshln}
\usepackage{youngtab}
\usepackage{young}
\usepackage{geometry}
\usepackage{bbm}
\usepackage{longtable}
\usepackage{url}
\usepackage{graphicx}
\usepackage{float}
\usepackage{units}
\usepackage{booktabs}
\usepackage{cite}
\usepackage{lineno}
\usepackage{subfigure}
\usepackage{indentfirst}
\usepackage{footnote}
\graphicspath{ {./images/} }
\usepackage[colorlinks,
linkcolor=blue,
anchorcolor=blue,
citecolor=blue,
urlcolor=black,
]{hyperref}
\usepackage{algorithm}
\usepackage{algorithmic}

\usepackage[]{caption2}

\renewcommand{\thesubfigure}{(\roman{subfigure})}
\makeatletter \renewcommand{\@thesubfigure}{\thesubfigure \space}
\renewcommand{\p@subfigure}{} \makeatother
\DeclareRobustCommand{\&}{%
    \ifdim\fontdimen1\font>0pt
    \textsl{\symbol{`\&}}%
    \else
    \symbol{`\&}%
    \fi
}
\newtheorem{theorem}{Theorem}[section]

\newtheorem{lemma}[theorem]{Lemma}

\newtheorem{proposition}[theorem]{Proposition}

\nolinenumbers
\title{A generalized Nystr\"{o}m method with column sketching for low-rank approximation of nonsymmetric matrices }
\author{Yatian Wang \footnotemark[2] \footnotemark[4] \and Hua Xiang\footnotemark[2] \footnotemark[3] \footnotemark[4] \and Chi Zhang\footnotemark[2] \and Songling Zhang\footnotemark[2]}

\begin{document}
    \nolinenumbers
    \date{}
    \maketitle
    \renewcommand{\thefootnote}{\fnsymbol{footnote}}
    \footnotetext[2]{School of Mathematics and Statistics, Wuhan University, Wuhan 430072, China.}
    \footnotetext[3]{Hubei Key Laboratory of Computational Science, Wuhan University, Wuhan 430072, China.}
    \footnotetext[4]{Corresponding author. E-mail address: {\tt yatianw@whu.edu.cn, hxiang@whu.edu.cn}.}

    \begin{abstract}
        This paper is concerned with the low-rank approximation for large-scale nonsymmetric matrices. Inspired by the classical Nystr\"{o}m method, which is a popular method to find the low-rank approximation for symmetric positive semidefinite matrices, we explore an extension of the Nystr\"{o}m method to approximate nonsymmetric matrices. The proposed method is a generalized Nystr\"{o}m method with column sketching and shows its advantages in accuracy and speed without sacrificing stability. And the numerical experiments will illustrate the robustness of our new methods in finding a desired low-rank approximation of nonsymmetric matrix.
    \end{abstract}

    \section{Introduction}

    As one core problem in numerical linear algebra, finding a low-rank approximation for a large-scale matrix $A\in\mathbb{R}^{m\times n}$ emerges in so many theoretical studies and practical cases \cite{markov2014,de2009}. There are a wide range of applications for low-rank approximation, including fluid dynamics \cite{dynamic2022}, astronomical imaging \cite{inference2021}, uncertainty quantification \cite{PG2015}, and other significant areas. Low-rank matrix approximations, such as the rank-revealing factorization, the truncated singular value decomposition (SVD), and the Nystr\"{o}m method, play a fundamental role in computational sciences and data analysis \cite{clarkson2009}. Among them, the Nystr\"{o}m method \cite{williams2000, nystrom1930} is highly successful for the low-rank approximation of symmetric positive semidefinite (SPSD) matrices such as Laplacian and kernel matrices that arise in machine learning applications \cite{gittens2013}. The Nystr\"{o}m method gives an accurate low-rank approximation for any symmetric matrix, when the singular value decaying speed is sufficiently fast. As for the low-rank approximation of general nonsymmetric matrices, the randomized singular value decomposition (rSVD) \cite{HMT2011}, rank-revealing QR decomposition \cite{demmel2015, Tropp2020}, and generilized Nystr\"{o}m method (GN) \cite{Naka2020} are ubiquitous and can be served as the useful and efficient tools in scientific research \cite{Naka2023, Tropp2023}. They appear frequently as the matrices having low numerical rank \cite{UT2019}.

    We focus on a generalization or extention of the classical Nystr\"{o}m method and aim at finding more efficient methods for the approximation of nonsymmetric matrices. Therefore, the existing ideas for using the Nystr\"{o}m method shall be reviewed briefly.

    \textit{Nystr\"{o}m methods for symmetric matrices.} Firstly let $A \in \mathbb{R}^{n \times n}$ be an SPSD matrix and the positive integer $k$ is the target rank and $r$ is called the sketching size. The most common form of Nystr\"{o}m methods is taken by $A_{nys}^{(r)} = WC^{\dagger}W^{T}$, where $W = A\Omega \in \mathbb{R}^{n\times r}$ and $C^{\dagger}$ is the pseudo-inverse of  the core matrix $C = \Omega^{T}A\Omega \in \mathbb{R}^{r\times r} (k\leq r\ll n)$ and $\Omega \in \mathbb{R}^{n\times r}$ is a randomized sketching matrix \cite{nystrom2023}. Among several variants of Nystr\"{o}m methods, there are two different truncated Nystr\"{o}m methods to obtain a target rank-$k$ approximation of $A_{nys}^{(r)}$. One version is given by $[A_{nys}^{(r)} ]_{k} = [ WC^{\dagger}W^{T}_{k}$, which performs rank-truncation in the Nystr\"{o}m approximation $A_{nys}^{(r)}$ directly \cite{wang2019, SPD2017}. The other traditional method defined by $ A_{nys}^{(r,k)}  = W [ C ]_{k}^{\dagger}W^{T}$ is to truncate the core matrix $C$ to target rank-$k$, which is cheaper to compute but takes less advantage of the full Nystr\"{o}m approximaton \cite{gittens2013, li2014}. It should be noted that $\lVert A - [ A_{nys}^{(r)} ]_{k} \lVert \leq \lVert A - A_{nys}^{(r,k)} \lVert $ holds in the Frobenius norm and the spectral norm for SPSD matrix \cite{pour2018}. As for an indefinite matrix, and the main concerned problem is the inaccuracy and instability when computing the pseudo-inverse of core matrix $C$, because the positive and negative eigenvalues of $A$ may \textquotedblleft cancel\textquotedblright\,each other out while forming $C$, which makes the eigenvalues of $C$ much smaller than $\sigma_{k}(A)$ by Cauchy's interlacing theorem \cite{Naka2023}. Additionally, the computation of $C^{\dagger}$ can be numerically unstable if $\sigma_{min}(C) < u$ where $u$ is the unit roundoff. In \cite{cai2022}, an error bound for the indefinite Nystr\"{o}m method is derived, which depends on how close the function values of the sampled points are. With the aim at improving the stability, a method that truncates the core matrix $\Omega^{T}A\Omega$ so that $\sigma_{min}((\Omega^{T}A\Omega)_{\epsilon}) > \epsilon$ where $\epsilon$ is of the order of the unit roundoff is proposed \cite{Naka2020}. Another version in \cite{gisb2015, oglic2018} makes use of the eigenvalue information of the original matrix, which is expensive to compute. The submatrix-shifted Nystr\"{o}m is suggested to deal with the indefinite matrix of few negative eigenvalues in \cite{ray2022}. Different from the methods described above using column sampling matrices for $\Omega$ \cite{wang2016}, the authors in \cite{Naka2023} utilize the randomized embeddings to overcome the instability. They also establish relative-error nuclear norm bounds of the resulting approximation for the case where the singular values decay rapidly \cite{Tropp2020,Naka2023}.

    \textit{Nonsymmetric approaches.} As for nonsymmetric and rectangular matrix $A \in\mathbb{R}^{m \times n}$, the rSVD is robust to approximate $A$ by the form $QQ^{T}A$, where the basis matrix $Q\in\mathbb{R}^{m\times r}$ is expected to contain as few orthonormal columns as possible and its range can approximate the range of $A$ \cite{HMT2011}. Nevertheless, the cost of computing $Q^{T}A$ (assuming $A$ is dense) becomes prohibitive when $n$ and $r$ are large. In \cite{demmel2023}, a generalized LU-factorization (GLU) for low-rank approximation is proposed and its good performance has been reported, as well. Meanwhile, generalizing the Nystr\"{o}m method to low-rank approximation for unsymmetric matrix is also significant and interesting, and it is natural to expect that more than one randomized matrix are needed to match and sketch from the right and left subspaces of $A$. In \cite{Naka2020}, a fast and near-optimal method called generalized Nystr\"{o}m method (GN) that guarantees the numerical stability is derived despite the presence of ill-conditioned pseudo-inverse. Meanwhile it leads to an approximation of $A$ with the form $AX(Y^{T}AX)^{\dagger}Y^{T}A$, where $X \in \mathbb{R}^{n \times r}$ and $Y \in \mathbb{R}^{m \times (r+l)}$ are two independent randomized sketching matrices ($l$ is an oversampling parameter). Moreover there is an inexpensive modification of GN to guarantee stability. The technique is to denote the $\epsilon$-pseudo-inverse of the core matrix as $(Y^{T}AX)^{\dagger}_{\epsilon} = V_{1}\Sigma^{-1}_{1}U^{T}_{1}$, where $Y^{T}AX = \begin{bmatrix}
        U_{1}, U_{2}
    \end{bmatrix}
    \begin{bmatrix}
        \begin{smallmatrix}
            \Sigma_{1} & 0 \\
            0 & \Sigma_{2}
        \end{smallmatrix}
    \end{bmatrix}
    \begin{bmatrix}
        V_{1}, V_{2}
    \end{bmatrix}^{T}$ is the SVD, $\Sigma_{1}$ contains the singular values larger than $\epsilon = \mathcal{O}(u\lVert A \lVert)$, and $U_{1}$ and $V_{1}$ consist of the left and right singular vectors corresponding to the singular values of $\Sigma_{1}$ \cite{Naka2020}. Another natural different version to form one small core matrix is $C = (AX)^{\dagger}A(Y^{T}A)^{\dagger}$, which minimizes the Frobenius norm of the error at the cost of more computational requirements \cite{cort2020}.

    Obviously when $A \succeq 0$ and $X = Y = \Omega$, the approximation form of GN is transformed into the classical Nystr\"{o}m method. Furthermore, the GN method reviewed above mainly concentrates on the decomposition of core matrix $C$ or the choice of the randomized sketching matrix. One shortcoming is that $A$ needs to be multiplied with two independent sketching matrices $X$ and $Y$ so that the numerical accuracy and stability may be affected.

    \textit{Related tools.} One of the core tools, which is widely utilized in random algorithms and has to be mentioned in randomized linear algebra, is randomized sketching \cite{FEVM}. Traditionally, the randomized matrices are chosen to be the column sampling matrices, which have exactly one nonzero entry equal to 1 in each column \cite{williams2000}. In this case, $AX$ or $Y^{T}A$ is a subset of $A$ and the core matrix $C$ is a principal submatrix of $A$. The common schemes for column sampling include uniform sampling \cite{wood2014}, k-means$++$ sampling \cite{oglic2017} and leverage score sampling \cite{musco2017, mahoney2009}, which are used in many scientific subjects and excellent works in recent years \cite{KLDP}. The application of randomized embedding is often referred to as sketching \cite{Tropp2020}, because the results obtained via randomized embeddings usually have smaller variance than the results yielded by column sampling.
    The common sketching methods include random partial isometry and Gaussian embedding, the latter being the simplest construction. In addition to the techniques mentioned above, the structured random embeddings are more practical and convenient to construct, to store, and to apply to vectors, such as sparse sign matrices \cite{wood2014}, subsampled trigonometric transforms \cite{bout2013, tropp2011}, and tensor random projections \cite{nelson2013}. There's no doubt that randomized embedding owns more attractive properties than column sampling matrices in many cases.

    \textit{Notations.} Throughout, we use dagger \textquotedblleft $\dagger$\textquotedblright\, to denote the pseudo-inverse of a matrix, and use $\lVert \cdotp\lVert_{F}$ for the Frobenius norm and $\lVert \cdotp\lVert_{2}$ for the spectral norm. Unless specified otherwise, denote $\sigma_{i}(A)$ as the $i$th largest singular value of $A$. Additionally, we use $\hat{A}_{r}$ to denote the rank-$r$ approximation of nonsymmetric matrices via GN or other approaches. We write \textquotedblleft $[Q,R] = \mathrm{qr\_econ}(A)$\textquotedblright\,to denote the \textit{economy-size} QR factorization of $A$ that only computes the first $n$ columns of $Q$ and the first $n$ rows of $R$, when $m \geq n.$ Arguments that are not needed are replaced by \textquotedblleft $\sim$\textquotedblright, so that, for instance, $[Q,\sim] = \mathrm{qr\_econ}(A)$ only returns the matrix $Q$. The \textquotedblleft iid \textquotedblright\,means independent and identically distributed.

    The remainder of this paper is organized as follows. In section 2, we briefly review the useful randomized techniques and the popular methods to obtain a low-rank approximation of nonsymmetric matrix, such as the generalized Nystr\"{o}m method and rSVD. The description and basic illustration of the algorithms proposed by us are shown in section 3. In addition, section 4 presents some theoretical analysis of the error bounds to discuss the accuracy of related methods. Numerical experiments to compare the performances and differences of GN, rSVD, and the new methods are shown in section 5, which also verify the error bound analysis. Finally, some conclusions will be drawn in section 6.

    \section{Generalized Nystr\"{o}m method and related techniques}

    In this section, firstly, we give a brief review of some related techniques, including the rank-revealing factorization and the randomized embedding, which are widely used in randomized linear algebra. For low-rank approximation of $A\in\mathbb{R}^{m\times n}$, these fundamental techniques are efficient, when the numerical target rank $k$ satisfies $k \ll min(m, n)$. Nextly, the existing methods, especially the GN and rSVD for low-rank approximation of nonsymmetric matrices will be introduced in detail.

    \subsection{Rank-revealing factorizations}
    At first, we will introduce the basic form of the rank-revealing factorization and show the procedures of the RURV algorithm as a practical approach to get a randomized rank-revealing factorization of $A$. Then another remarkable technique, the randomized embedding will be described as well.

    The rank-revealing factorizations have been used to solve ill-conditioned linear systems and least-squares problem. Firstly, let us investigate what it means when denoting a factorization as \textit{rank-revealing}. Given a matrix $A \in \mathbb{R}^{m \times n}$, its factorization is provided with the form
    \begin{equation}
        A = URV^{T},
    \end{equation}
    where $U \in \mathbb{R}^{m \times c}$, $R \in \mathbb{R}^{c\times n}$, $V \in \mathbb{R}^{n \times n}$, and $c = \mathrm{min}(m, n)$. The matrices $U$ and $V$ are orthogonal and $R$ is upper-triangular (or banded upper-triangular). We expect this factorization to reveal the numerical rank of $A$ in the sense of that we can obtain a near-optimal approximation of $A$ by truncating (1) to any level $k$. That is,
    \begin{equation}
        \lVert A - U(: , 1 : k)R(1 : k, :)V^{T} \lVert \approx \inf{ \lVert A - B \lVert},
    \end{equation}
   where $B$ has rank $k$, for $k = {1, 2, \cdots , c}$ \cite{Tropp2020}.

   Obviously, the classical deterministic techniques for computing the factorization are usually expensive. One considerable algorithm proposed by Demmel, Dumitriu and Holtz is efficient to compute a rank-revealing URV factorization of an $m \times n$ matrix $A$, typically with $m \geq n$.  This method is called RURV and it proceeds as follows \cite{demmel2007}.
    \begin{itemize}
         \raggedright
         \item [(\romannumeral1)]
         Draw a randomized test matrix $\Omega$ from a standard normal distribution;
         \item [(\romannumeral2)]
         Perform an unpivoted QR factorization of $\Omega$: $[V, \sim] = \mathrm{qr\_econ}(\Omega)$;
         \item [(\romannumeral3)]
         Perform an unpivoted QR factorization of $AV$: $[U, R] = \mathrm{qr\_econ}(AV)$.
    \end{itemize}

   The steps (\romannumeral1) and (\romannumeral2) generate the matrix $V$ whose columns can serve as a random orthonormal basis. It is easy to verify that the matrices $U$, $R$, and $V$ satisfy the factorization in (1), the cost of which is dominated by two unpivoted QR factorizations \cite{Ballard2019}. To improve the rank-revealing ability of the factorization, we can incorporate a small number of power iteration in the step 2: $[V, \sim] = \mathrm{qr}((A^{T}A)^{q}\Omega)$, where $q = 1$ or $q = 2$ in practice \cite{gopal2018}. In \cite{HMT2011}, the technique of constructing an orthonormal matrix $Q$ whose range approximates the range of $A$ is actually consistent with the randomized rank-revealing factorization.

    \subsection{Randomized embedding}
    The performance of the randomized embedding is attractive, which preserves the 2-norm of every vector in a given subspace. Given $A \in \mathbb{R}^{m \times n}$, $x \in \mathbb{R}^{n}$, and a distortion parameter $\epsilon \in (0, 1)$, then a linear map $S \in \mathbb{R}^{s\times m}$, which enacts a dimension reduction is a subspace embedding with distortion $\epsilon$ when
    $$(1 - \epsilon)\lVert Ax \lVert \leq \lVert SAx \lVert \leq (1+ \epsilon)\lVert Ax \lVert,$$
    where $s\ll m$, which means that the embedding $S$ transfers data from the high-dimensional space to the low-dimensional space \cite{wood2014}. Morever, the randomized emdedding has more attractive properties than column sampling matrices \cite{Tropp2020}. As the most widely used randomized embedding, Gaussian embedding has optimal theoretical guarantees and is frequently used in theoretical analysis.

    \noindent\textbf{Gaussian matrices.} A Gaussian embedding is a random sketching matrix of the form
    $$G \in \mathbb{R}^{s \times m},\quad (G)_{ij} \sim \mathcal{N}(0,1/s).$$
    For each $x \in \mathbb{R}^{m}$, the scaling of the matrix ensures that $\mathbb{E}[\lVert Gx \lVert^{2}] = \lVert x \lVert^{2}$. The cost of applying a random Gaussian matrix to an $m \times n$ matrix is $\mathcal{O}(smn)$. Meanwhile, the Gaussian embedding admits a simple theoretical analysis, which tends to be difficult for other randomized embeddings though they exhibit similar behavior to the Gaussian embeding in practical experiences. So for many practical purposes, Gaussian analysis can provide promising insight and is usually used to provide a rule of thumb for the general behavior \cite{Tropp2020, Naka2023}.

     \noindent\textbf{Sparse maps.} Next, we describe the sparse dimension reduction map, whose entries are random signs \cite{nelson2013, HMT2011}. The sparse sign matrix $S \in \mathbb{R}^{s \times m}$ with nonzero entries is useful for sparse data and may require less data movement. An effective construction of it with a sparsity parameter $\zeta$ takes the form
     $$ S = \sqrt{\frac{m}{\zeta}} [s_{1}, \cdots , s_{m}] \in \mathbb{R}^{s \times m},$$
     where each column $s_{i}$ of $S$ is iid random vector, which is constructed by drawing $\zeta$ iid random signs and has exactly $\zeta$ nonzero entries that take $\pm1$ with equal probability, placed uniformly at random coordinates \cite{Naka2024}. The sketching size $s = \mathcal{O}(r\log r)$ is chosen for theoretical guarantees in \cite{Naka2023, cohen2016}. In \cite{tropp2019}, the sparsity level is recommended as $\zeta = \min\left\lbrace s, 8\right\rbrace$ in practice.

     \noindent\textbf{SRTTs.} The subsampled randomized trigonometric transform matrix takes the form
     $$S = \sqrt{\frac{m}{s}}DFR^{T} \in \mathbb{R}^{m\times s}.$$
     In this expression, $D \in \mathbb{R}^{m \times m}$ is a diagonal matrix, whose entries are chosen independently at random and take $\pm1$ with equal probability. $F \in \mathbb{R}^{m \times m}$ is commonly the discrete cosine transform (DCT) in the real case \cite{Naka2023}. And $R \in \mathbb{R}^{s \times m}$ is a random restriction. The sketching size needs to be $s = \mathcal{O}(r\log r)$ \cite{tropp2011}. The cost of applying SRTT to an $m \times n$ matrix is $\mathcal{O}(mn\log r)$ by the subsampled FFT algorithm in \cite{FEVM}.

    \subsection{Generalized Nystr\"{o}m method}
    As mentioned above, to obtain a rank-$r$ approximation $\hat{A}_{r}$ of a nonsymmetric matrix $A \in \mathbb{R}^{m \times n}$, a natural starting point is to look for an approximation with column space $AX$ and row space $Y^{T}A$ for random sketching matrices $X \in \mathbb{R}^{n \times r}$ and $Y \in \mathbb{R}^{m \times r}$. We refer to the specific procudures in Algorithm 2-1 as the generalized Nystr\"{o}m method, or GN for short. The rank-$r$ approximation of $A$ by GN takes the form
        \begin{equation}
            A \approx \hat{A}_{GN} = AX(Y^{T}AX)^{\dagger}Y^{T}A.
        \end{equation}

    This stable and fast method has near-optimal approximation quality and the computational cost is $\mathcal{O}(mn\log n +r^{3})$ for dense matrices. Then let us introduce the pseudocodes of GN and discuss its implementation details.
    \begin{algorithm}[htb]
        \renewcommand{\thealgorithm}{2-1}
        \renewcommand{\algorithmiccomment}[1]{\hfill $\triangleright$ #1}
        \caption{GN \cite{Naka2020}: generalized Nystr\"{o}m approximation for nonsymmetric matrices.}
        \begin{algorithmic}[1]
            \REQUIRE
            The nonsymmetric matrix $A \in \mathbb{R}^{m \times n}$, $r\in\mathbb{N}$, oversampling $l = \lceil 0.5r \rceil$
            \ENSURE
            A rank-$r$ approximation $\hat{A}_{GN}$
            \STATE Draw the sketching matrices $X \in \mathbb{R}^{n \times r}$, $Y \in \mathbb{R}^{m \times (r+l)}$;
            \COMMENT{Gaussian embedding}
            \STATE Compute $AX$, $Y^{T}A$, $Y^{T}AX$;\\
            \STATE $[Q, R]$ = qr\_econ$(Y^{T}AX)$;\\
            \STATE $\hat{A}_{GN} = ((AX)R^{-1})(Q^{T}(Y^{T}A))$; \COMMENT{GN}\\
        \end{algorithmic}
    \end{algorithm}

    The GN performs a straightforward QR factorization of the core matrix. Nevertheless, another generalized Nystr\"{o}m method for rank-$r$ approximation of $A$ with an inexpensive modification is called stabilized generalized Nystr\"{o}m, or stabilized GN \cite{Naka2020}. The way to deal with the core matrix $Y^{T}AX$ of stabilized GN is different from that of GN.
    The stabilized GN denotes a $\varepsilon$-pseudo-inverse of the core matrix $Y^{T}AX$, whose form is $(Y^{T}AX)^{\dagger}_{\epsilon} \approx V_{1}\Sigma^{-1}_{1}U^{T}_{1}$, where $Y^{T}AX = \begin{bmatrix}
        U_{1}, U_{2}
    \end{bmatrix}
    \begin{bmatrix}
        \begin{smallmatrix}
            \Sigma_{1} & 0 \\
            0 & \Sigma_{2}
        \end{smallmatrix}
    \end{bmatrix}
    \begin{bmatrix}
        V_{1}, V_{2}
    \end{bmatrix}^{T}$ is the SVD,
    $\Sigma_{1}$ contains singular values larger than $\epsilon = \mathcal{O}(u\lVert A \lVert)$ and the $U_{1}$ and $V_{1}$ consist of the left and right singular vectors respectively, which are corresponding with the singular values contained in $\Sigma_{1}$. Another cheaper alternative with the core matrix is to compute the QR factorization $[Q,R] = qr(Y^{T}AX)$ and look for the diagonal elements of $R$ less than $\epsilon$ \cite{Naka2020}.

    From the procudures of Algorithm 2-1, the cost to calculate the $AX(R)^{-1}$ or $Q^{T}(Y^{T}A)$ is $\mathcal{O}(mr^{2})$ and the memory requirement is $\mathcal{O}((m + 1.5n)r)$. And it needs to be noted that the output of GN fails to give an approximative truncated SVD because no orthogonal columns are obtained. Nonetheless, the GN method shows its advantages in stability and faster speed in computation. Then we present the procedures of rSVD in Algorithm 2-2 and denote the rank-$r$ approximation via rSVD as $\hat{A}_{rSVD}$.
    \begin{algorithm}[H]
        \renewcommand{\thealgorithm}{2-2}
        \renewcommand{\algorithmiccomment}[1]{\hfill $\triangleright$ #1}
        \caption{ rSVD \cite{HMT2011}: randomized SVD approximation for nonsymmetric matrices.}
        \begin{algorithmic}[1]
            \REQUIRE
            The nonsymmetric matrix $A \in \mathbb{R}^{m \times n}$, $r\in\mathbb{N}$
            \ENSURE
            A rank-$r$ approximation $\hat{A}_{rSVD}$
            \STATE Draw the sketching matrices $\Omega \in \mathbb{R}^{n \times r}$;
            \COMMENT{Gaussian embedding}
            \STATE Compute $A\Omega$;\\
            \STATE Orthogonalize $A\Omega$ to obtain $Q$: $[Q, \sim]$ = qr\_econ$(A\Omega)$;\\
            \STATE $[U_{0}, \Sigma_{0}, V_{0}]$ = svd$(Q^{T}A)$; \\
            \STATE $\hat{A}_{rSVD} = (QU_{0})\Sigma_{0}V_{0}^{T}$.  \COMMENT{rSVD}\\
        \end{algorithmic}
    \end{algorithm}
    The cost of Algorithm 2-2 is generally dominated by the product $Q^{T}A$ and $A\Omega$, which require $\mathcal{O}(rmn)$ flops for dense matrices. And the steps 4 and 5 take $\mathcal{O}(r^{2}n)$ and $\mathcal{O}(r^{2}m)$ flops, respectively.

     In order to illustrate the differences between these methods adequately, the GN and rSVD will be performed in the following numerical experiments, and the error bound analysis for them will be shown in section 4 as well. In the following part of this section, it is still necessary to compare GN and rSVD with other existing algorithms and summarize the main advantages and disadvantages of them briefly.

     In \cite{demmel2023}, from the perspective of rank-revealing LU factorization, the authors derive an approximation that takes the form $A \approx TS$, where $T$ is a tall-skinny matrix and $S$ is a short-wide matrix. The numerical experiments suggest that its accuracy is between that of GN and rSVD, and the arithmetric cost is more expensive than GN. In \cite{gu2015}, the version of subspace iteration method is used and it is obvious that the flop cost of subspace iteration method is much higher than GN. The methods in \cite{Tropp2017} and \cite{clarkson2009} are mathematically the same, but the approach in \cite{Tropp2017} uses orthogonalization to make sure a stable output. And the approximation form in \cite{clarkson2009} is essentially equivalent to that of GN. So after careful analysis, the GN and rSVD possess more advantages in approximating nonsymmetric matrices so that we will mainly compare our methods with them. Certainly, other methods mentioned above behave successfully in so many applications and play an important role in finding randomized low-rank matrix approximation as well \cite{markov2014}.

    \section{Our proposed algorithms}\setcounter{footnote}{0}\renewcommand{\thefootnote}{\arabic{footnote}}
    In light of the observations and analysis of these earlier-introduced methods, our primary idea for approximating a nonsymmetric matrix is to discover other efficient approaches to deal with the core matrix with more stability and make the output $\hat{A}_{r}$ access to information available with orthonormal columns. Though the GN avoids the presence of an ill-conditioned pseudo-inverse, but the inverse of $R$ may cause instability for large-sparse matrices. And we can merely try to compute the QR decomposition of $AX$ and $Y^{T}A$ and calculate the SVD of $R_{1}(Y^{T}AX)^{\dagger}R_{2}^{T}$, where $R_{1}$ and $R_{2}$ are the R-factors of $AX$ and $A^{T}Y$, which can generate an approximate SVD of $A$. Obviously, the computational cost of this method is extremely expensive and do not exhibit enough advantages in comparison with GN and rSVD. Take these difficulties into consideration, it is necessary to make a difference in the way of row-sketching and cloumn-sketching and
    discover other efficient and innovative approaches to deal with the core matrix.
    \subsection{A primary algorithm GN-r\&c}
    Inspired by the RURV algorithm \cite{demmel2007}, which is a simple randomized algorithm for computing a rank-revealing factorization of a general matrix, whose detailed procedures have been introduced in section 2. And the step (\romannumeral2) of RURV algorithm that performs QR decomposition of one Gaussian matrix, generates a matrix $V$, whose columns can serve as a randomized orthonormal basis. Therefore, it motivates us to look for the corresponding randomized orthonormal matrices, which can be employed to update the initial Gaussian matrices $X$ and $Y$ in the core matrix.

    Then for $A\in\mathbb{R}^{m\times n}$, it is natural to take $AX$ and $A^{T}Y$ as the column space and row space of $A$, respectively. And the orthonormal matrices $Q_{1}$ and $Q_{2}$ are the Q-factors of $AX$ and $A^{T}Y$, whose ranges can approximate the ranges of $A$ and $A^{T}$, respectively. So it is reliable to update $X$ and $Y$ with $Q_{2}$ and $Q_{1}$. Thus we have
    \begin{equation}
        A \approx (AQ_{2})(Q_{1}^{T}AQ_{2})^{\dagger}(A^{T}Q_{1})^{T}.
    \end{equation}

    It is reasonable to make such modifications, which have the potential to make the approximation more accurate. And if we continue to perform the QR factorization of $AQ_{2}$, that is, $AQ_{2} = \tilde{Q}\tilde{R}$, we can generate the pseudo-inverse of the core matrix with the form
    $$(Q_{1}^{T}AQ_{2})^{\dagger} = \tilde{R}^{\dagger}\tilde{Q}^{T}Q_{1},$$
    which brings further computational stability.
    Now we show the pseudocodes of this improved method in Algorithm 3-1, and we refer to it as generalized Nystr\"{o}m method with row and column sketching, or GN-r\&c. We denote the rank-$r$ approximation of nonsymmetric matrices by GN-r\&c as $\hat{A}_{GN\mbox{-}r\&c}$
    \begin{algorithm}[htb]
        \renewcommand{\thealgorithm}{3-1}
        \renewcommand{\algorithmiccomment}[1]{\hfill $\triangleright$ #1}
        \caption{GN-r\&c: generalized Nystr\"{o}m approximation for nonsymmetric matrices with both row and column sketching.}
        \begin{algorithmic}[1]
            \REQUIRE
            The nonsymmetric matrix $A \in \mathbb{R}^{m \times n}$, $r\in\mathbb{N}$
            \ENSURE
            A rank-$r$ approximation $\hat{A}_{GN-r\&c}$
            \STATE Draw the sketching matrices $X \in \mathbb{R}^{n \times r}$, $Y \in \mathbb{R}^{m \times r}$;
            \COMMENT{Gaussian embedding}
            \STATE Compute $AX$, $A^{T}Y$;\\
            \STATE $[Q_{1}, \sim]$ = qr\_econ$(AX)$, $[Q_{2}, \sim]$ = qr\_econ$(A^{T}Y)$;
            \COMMENT{generate random orthonormal matrices}
            \STATE Compute $AQ_{2}$, $A^{T}Q_{1}$;
            \COMMENT{update $X$ and $Y$ with $Q_{2}$ and $Q_{1}$}\\
            \STATE $[\tilde{Q}, \tilde{R}]$ = qr\_econ$(AQ_{2})$, or $[\hat{Q}, \hat{R}]$ = qr\_econ$(A^{T}Q_{1})$;\\
            \STATE $\hat{A}_{GN-r\&c} = (AQ_{2})(Q_{1}^{T}\tilde{Q}\tilde{R})^{\dagger}(A^{T}Q_{1})^{T}$, or $(AQ_{2})(\hat{R}^{T}\hat{Q}^{T}Q_{2})^{\dagger}(A^{T}Q_{1})^{T}$;\COMMENT{GN-r\&c}\\
        \end{algorithmic}
    \end{algorithm}

    Analogous to the steps 1-2 in Algorithm 2-1, GN-r\&c generates two independent Gaussian matrices and the GN-r\&c in step 3 performs two QR factorizations of $AX$ and $A^{T}Y$. The improvement lies in making full use of the Q-factors of $AX$ and $A^{T}Y$ to update the original sketching matrices $Y$ and $X$, respectively. For steps 5 and 6, there are some remarks worthy to illustrate.
     \begin{itemize}
        \item[(\romannumeral1)]
        $Q_{1} \in \mathbb{R}^{m \times r}$ and $Q_{2} \in \mathbb{R}^{n \times r}$ are two independent random orthogonal matrices, whose ranges align well with the approximative ranges of $A$ and $A^{T}$, respectively. Indeed, step 5 contains another two separate randomized rank-revealing factorizations $\tilde{Q}\tilde{R}Q_{2}^{T}$ of $A$ and $\hat{Q}\hat{R}Q_{1}^{T}$ of $A^{T}$.
        \item[(\romannumeral2)]
        Unlike the way to compute the SVD or QR factorization of the initial core matrix $Y^{T}AX$ directly, for GN-r\&c, $Q_{1}^{T}AQ_{2} = (Q_{1}^{T}\tilde{Q})(\tilde{R})$, or $(\hat{R}^{T})(\hat{Q}^{T}Q_{2})$, and the pseudo-inverse of which is $\tilde{R}^{\dagger}(\tilde{Q}^{T}Q_{1}$), or $(Q_{2}^{T}\hat{Q})(\hat{R}^{T})^{\dagger}$. So in our numerical experiments, the rank-r approximation of $A$ is taken actually by the form\\
        $A \approx (AQ_{2})(Q_{1}^{T}AQ_{2})^{\dagger}(A^{T}Q_{1})^{T} = (AQ_{2})\tilde{R}^{\dagger}(\tilde{Q}^{T}Q_{1})(A^{T}Q_{1})^{T}$,
        or $(AQ_{2})(Q_{2}^{T}\hat{Q})(\hat{R}^{T})^{\dagger}(A^{T}Q_{1})^{T}$, which avoids computing the pseudo-inverse of the core matrix that may be ill-conditioned.
        \item [(\romannumeral3)]
        We can also obtain two different approximation forms with orthonormal columns. It is easy to derive that
        \[A \approx (AQ_{2})(Q_{1}^{T}AQ_{2})^{\dagger}(A^{T}Q_{1})^{T} = \tilde{Q}\tilde{R}(Q_{1}^{T}\tilde{Q}\tilde{R})^{\dagger}\hat{R}^{T}\hat{Q}^{T} = (\tilde{Q}\tilde{R}\tilde{R}^{\dagger})(Q_{1}^{T}\tilde{Q})^{T}(\hat{R}^{T}\hat{Q}^{T}),\]
        \[A \approx (AQ_{2})(Q_{1}^{T}AQ_{2})^{\dagger}(A^{T}Q_{1})^{T} = \tilde{Q}\tilde{R}(\hat{R}^{T}\hat{Q}^{T}Q_{2})^{\dagger}\hat{R}^{T}\hat{Q}^{T} = \tilde{Q}\tilde{R}Q_{2}^{T}\hat{Q}((\hat{R}^{T})^{\dagger}\hat{R}^{T})\hat{Q}^{T}.\]
        The $\tilde{R}\tilde{R}^{\dagger}\tilde{Q}^{T}Q_{1}\hat{R}$ and $\tilde{R}Q_{2}^{T}\hat{Q}(\hat{R}^{T})^{\dagger}\hat{R}^{T}$ can be implemented via the  SVD to obtain the approximative SVD of $A$.
        \item [(\romannumeral4)]
        When assuming that $\tilde{R}$ is of full-row rank and $\hat{R}^{T}$ is of full-column rank, we have $\tilde{R}\tilde{R}^{\dagger} = I$ and $(\hat{R}^{T})^{\dagger}\hat{R}^{T} = I$, and then the approximation forms in (3) can be expressed as
        \[A \approx Q_{1}\hat{R}^{T}\hat{Q}^{T},\quad A \approx \tilde{Q}\tilde{R}Q_{2}^{T}.\]
        The simpler forms derived under hypothesis are coincident with rSVD.
        \item [(\romannumeral5)]
        The singular value decaying speed makes a remarkable affection on the convergency speed and accuracy of all the methods including GN-r\&c.
    \end{itemize}
     \begin{figure} [H]
        \subfigure[\textit{Cry10000}]{
            \includegraphics[width=8cm,height=6cm]{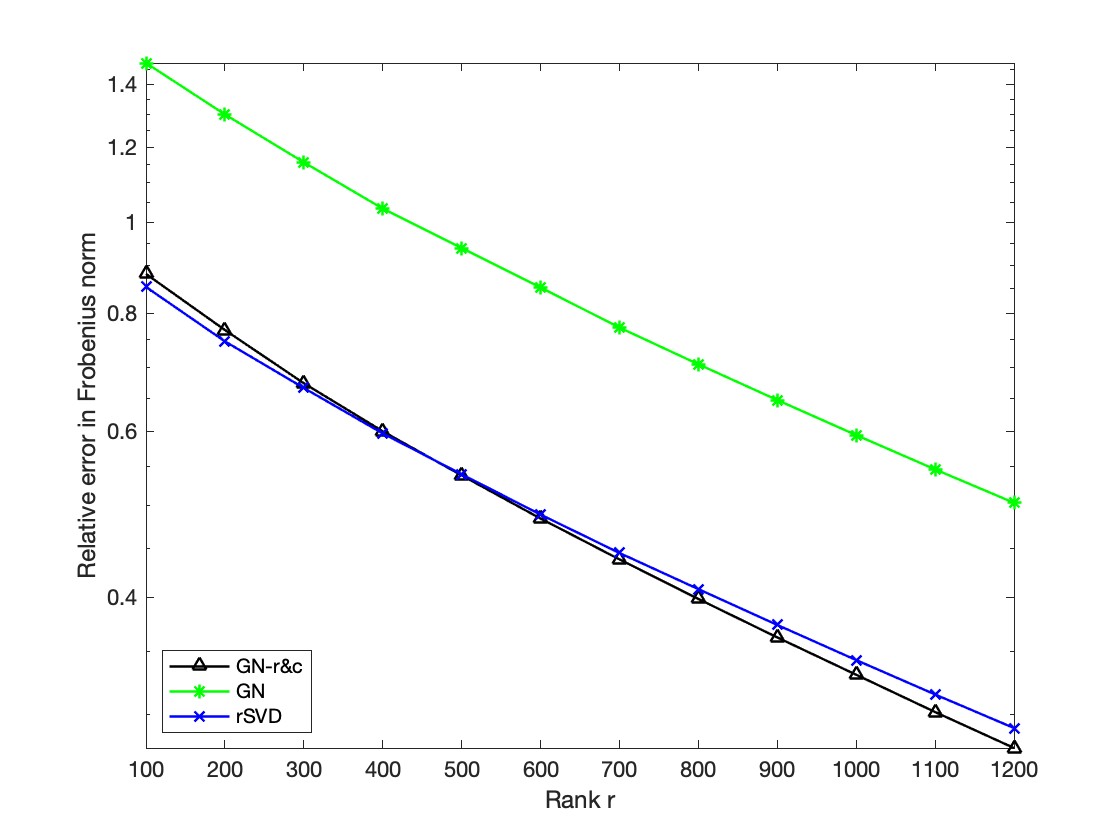} \label{Fig.sub.a2}
        }
        \hspace{2mm}
        \subfigure[\textit{Cry10000}]{
            \includegraphics[width=8cm,height=6cm]{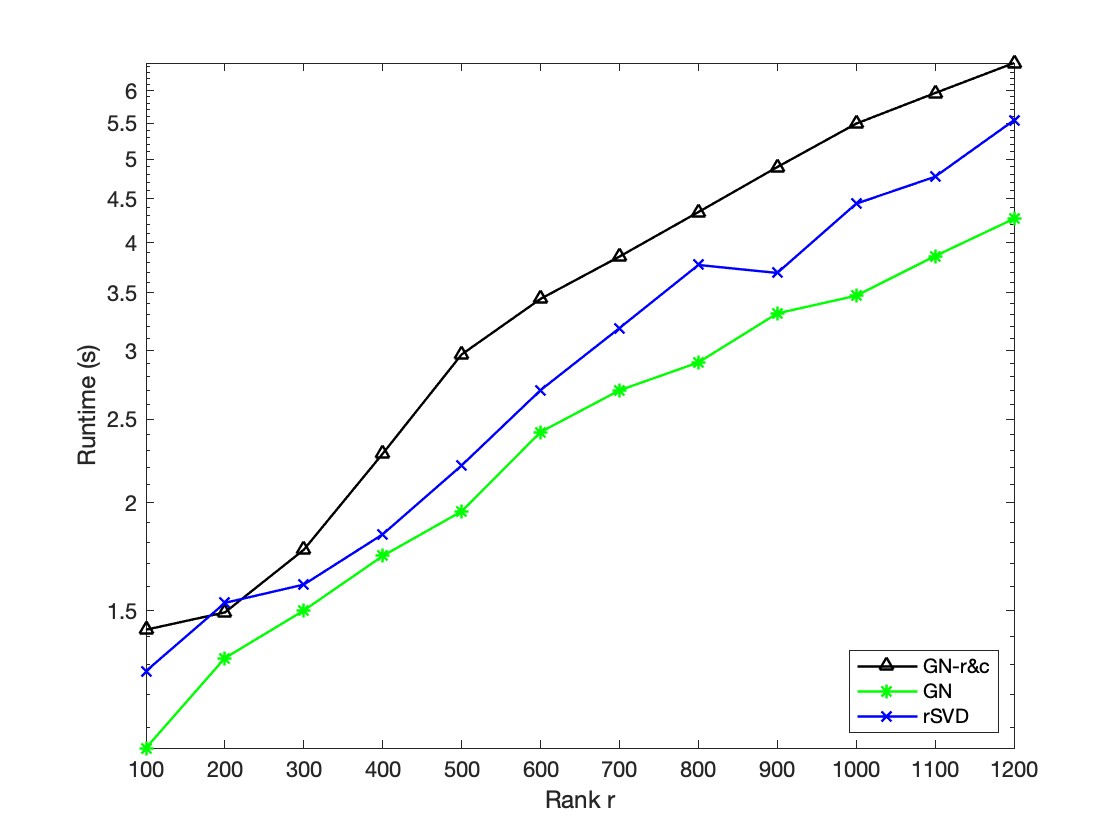} \label{Fig.sub.b2}
        }

        \subfigure[\textit{Synthetic matrix}]{
            \includegraphics[width=8cm,height=6cm]{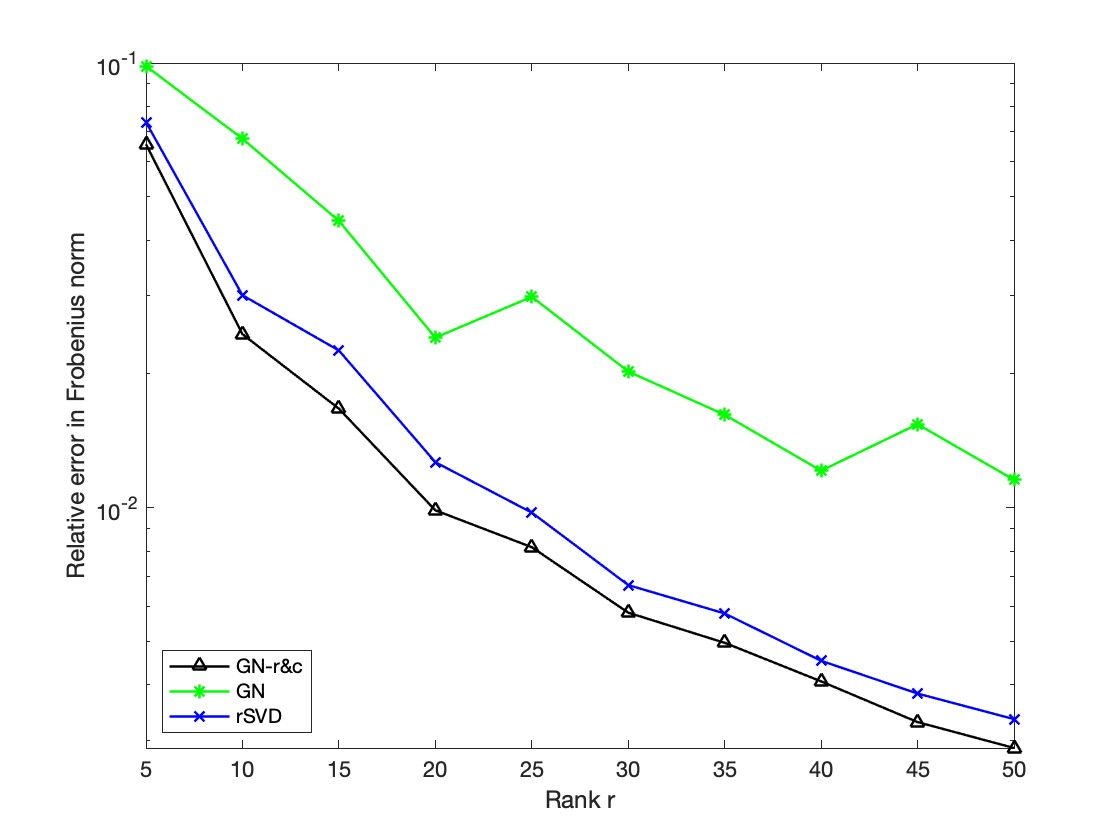} \label{Fig.sub.c2}
        }
        \hspace{2mm}
        \subfigure[\textit{Synthetic matrix}]{
            \includegraphics[width=8cm,height=6cm]{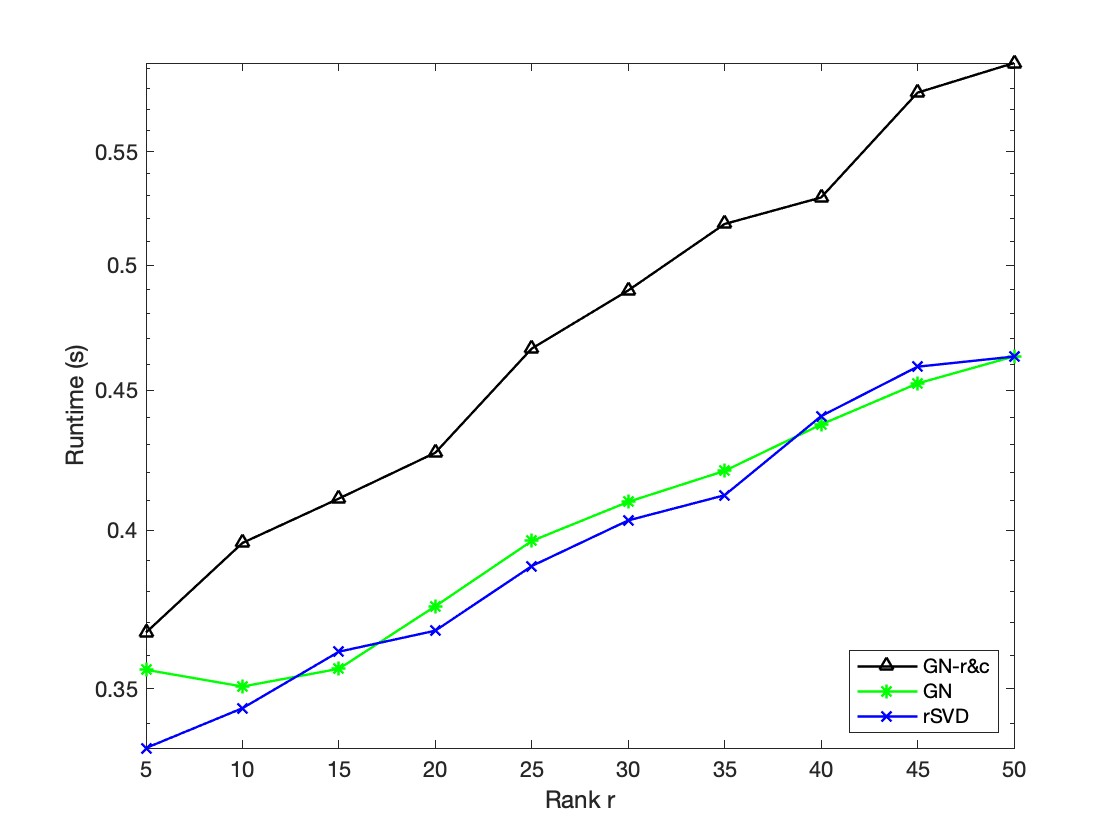} \label{Fig.sub.d2}
        }
        \caption{\textit{\textbf{Relative error and runtime of GN-r\&c, GN, and rSVD to approximate nonsymmetric matrices.} \ref{Fig.sub.a2} and \ref{Fig.sub.c2} show the relative error in Frobenius norm. The GN-r\&c performs better than GN and rSVD in stability and accuracy. In \ref{Fig.sub.b2} and \ref{Fig.sub.d2}, the speed of GN-r\&c is the slowest, especially for the dense matrix. The oversampling parameter of GN is 5 for the synthetic matrix and the recommended choice of $l$ for matrix Cry10000 is $0.5r$, because of different singular value decaying speeds.}}\label{Fig.main.2}
        \end {figure}

     From the observations of pseudocodes above, it is apparent that the cost of GN-r\&c is more expensive in comparison with the GN and rSVD. As for the accuracy, we try to choose the sparse and dense matrices separately to show the relative error $\frac{\lVert A - \hat{A}_{r}\lVert_{F}}{\lVert A\lVert_{F}}$ and the runtime of the GN, GN-r\&c, and rSVD, where the $\hat{A}_{r}$ is a rank-$r$ approximation obtained via GN or other methods.

    The nonsymmetric sparse matrix Cry10000 is $10000\times 10000$ and has 49699 nonzeros, whose singular spectrum is flat\footnote[1]{https://math.nist.gov/MatrixMarket/}. The other experiment is concerned with a $8000\times 8000$ synthetic dense matrix, which is generated by the form $A = M\Lambda N^{T}$, where $M$ and $N$ are all random orthogonal matrices (Q-factors of the random square Gaussian matrices) and $\Lambda$ has fast polynomial decaying singular values\footnote[2]{Throughout, all numerical experiments were performed in MATLAB (version 2022a) on a MacBook Pro with a 2.3GHz Intel Core i7 processor with four cores.}, where the singular value $\sigma_{i} = \frac{1}{i^{2}}$.

    A preliminary testing is performed in Fig.\ref{Fig.main.2}. It is not abrupt to find that GN-r\&c has a better approximation quality than GN and rSVD, even the speed of GN-r\&c is much slower. Meanwhile, the results indicate that the primary method GN-r\&c has made a steady progress in accuracy, but is not an extremely attractive algorithm due to the less ideal computational speed and much more storage requirement.
    \subsection{An improved algorithm GN-c}
    Encouraged by the good performance in accuracy and related analysis of the GN-r\&c in steps 5 and 6, we attempt to explore one faster, more accurate, but less complicated method that makes full use of the sketching matrix and its multiplications with $A$. From the procedures of Algorithm 3-1, GN-r\&c generates two independent sketching matrices and proceeds two randomized rank-revealing factorizations. Though GN-r\&c avoids calculating the pseudo-inverse of the core matrix $Q_{1}^{T}AQ_{2}$ directly, which is beneficial for enhanceing the stability, the bottleneck to speed up the calculation of GN-r\&c is four matrix-matrix multiplications and four involved QR factorizations.

    Furthermore, the two randomized rank-revealing factorizations $\tilde{Q}\tilde{R}Q_{2}^{T}$ of $A$ and $\hat{Q}\hat{R}Q_{1}^{T}$ of $A^{T}$ for GN-r\&c are separated from each other and help to derive two equivalent approximation forms. All of these facts lead us to consider that it can be enough to select merely one sketching matrix ($X$ or $Y$) and perform one randomized rank-revealing factorization for obtaining a low-rank approximation of $A$.

    For GN-r\&c, the ranges of the orthonormal matrices $Q_{1}$ and $Q_{2}$ approximate the ranges of $A$ and $A^{T}$ respectively. So simply for instance, if we only draw one Gaussian matrix $X$ in step 1 and update the \textquotedblleft $Y$\textquotedblright\,with the Q-factor $(Q_{1})$ of $AX$, we can obtain the form
    $$AX(Q_{1}^{T}AX)^{\dagger}(A^{T}Q_{1})^{T}.$$

     Along the same lines, update $X$ with the Q-factor of $A^{T}Q_{1}$ so that another more stable and accurate method arises. Up to now, owing to reasonable analysis and related observations above, it is feasible to generate a more efficient algorithm, called generalized Nystr\"{o}m method with column sketching, or GN-c. We present the whole procedures directly as follows.

     Firstly, we generate a randomized Gaussian matrix $X$, and perform the QR factorization of $AX$ in step 3 so that the columns of the $Q_{1}$ (Q-factor of $AX$) can serve as one sketching orthonormal basis. Then it is reasonable to substitute the previous Gaussian matrix $Y$ with the randomized orthonormal matrix $Q_{1}$. Indeed, the steps 3-5 can be seen as one whole randomized rank-revealing factorization $\hat{Q}\hat{R}Q_{1}^{T}$ of $A^{T}$. Lastly, we take the form $A\hat{Q}(\hat{R}^{T})^{\dagger}(A^{T}Q_{1})^{T}$ to gain a low-rank approximation of $A$. Meanwhile, we can start with the row sketching to approximate $A$, such as computing the $A^{T}Y$ in step 2, where $Y$ is a sketching matrix and the procedures are almost consistent with the GN-c method with column sketching.
     \begin{algorithm}[H]
        \renewcommand{\thealgorithm}{3-2}
        \renewcommand{\algorithmiccomment}[1]{\hfill $\triangleright$ #1}
        \caption{GN-c: generalized Nystr\"{o}m approximation for nonsymmetric matrices with column sketching.}
        \begin{algorithmic}[1]
            \REQUIRE
            The nonsymmetric matrix $A \in \mathbb{R}^{m \times n}$, $r\in\mathbb{N}$
            \ENSURE
            A rank-$r$ approximation $\hat{A}_{GN-c}$
            \STATE Draw the sketching matrices $X \in \mathbb{R}^{n \times r}$;
            \COMMENT{Gaussian embedding}
            \STATE Compute $AX$;\\
            \STATE Orthogonalize $AX$: $[Q_{1}, R_{1}]$ = qr\_econ$(AX)$;\\
            \STATE Compute $A^{T}Q_{1}$; \COMMENT{equal to the effect of $A^{T}Y$}\\
            \STATE Orthogonalize $A^{T}Q_{1}$: $[\hat{Q}, \hat{R}]$ = qr\_econ$(A^{T}Q_{1})$; \COMMENT{randomized rank-revealing factorization}\\
            \STATE $\hat{A}_{GN-c} = A\hat{Q}(Q_{1}^{T}A\hat{Q})^{\dagger}(A^{T}Q_{1})^{T} = A\hat{Q}(\hat{R}^{T})^{\dagger}(A^{T}Q_{1})^{T}$. \COMMENT{GN-c}\\
        \end{algorithmic}
    \end{algorithm}
    Necessarily, some remarks about GN-c and the connection with GN-r\&c should be emphasized as follows.
    \begin{itemize}
        \item[(\romannumeral1)]
        The GN-c method includes one initial Gaussian matrix and two QR factorizations. The steps 3-5 of GN-c derive a randomized rank-revealing factorization of $A^{T}$. In comparison with GN-r\&c, the computational cost of matrix-matrix multiplications and involved QR factorizations decrease in half. And the GN-c is more efficient to get a low-rank approximation, which will be verified in the numerical illustration.
        \item[(\romannumeral2)]
        For GN-c, the primary form of the core matrix $Y^{T}AX$ is simplified with updating the $X$ and $Y$ with $\hat{Q}$ and $Q_{1}$, and then we have
            \[Y^{T}AX = Q_{1}^{T}A\hat{Q} = \hat{R}^{T}.\]
        \item[(\romannumeral3)]
        If we do not update the $X$ with $\hat{Q}$ in step 6, we can generate
        \[A \approx AX(Q_{1}^{T}AX)^{\dagger}(A^{T}Q_{1})^{T} = Q_{1}R_{1}(R_{1})^{\dagger}(\hat{Q}\hat{R})^{T} = Q_{1}R_{1}(R_{1})^{\dagger}\hat{R}^{T}\hat{Q}^{T}.\]
        Assuming that $R_{1}$ is of full-row rank, which means $R_{1}(R_{1})^{\dagger} = I$, thus the equality above can be formed as
        \[A \approx Q_{1}\hat{R}^{T}\hat{Q}^{T},\] which is consistent with rSVD as well.
         \item[(\romannumeral4)]
          In this work, we are not aim to find a fixed rank-$r$ approximation of $A$. So GN-r\&c and GN-c do not need to consider an oversampling parameter in numerical tests.
    \end{itemize}

    The remarks (\romannumeral1) and (\romannumeral2) illustrate that GN-c is faster than GN-r\&c, and remarks (\romannumeral2) and (\romannumeral3) suggest that GN-c may be more accurate than rSVD and GN-r\&c. For GN-c, the Gaussian matrix is still chosen in the sketching technique because of its strong theoretical guarantees.

    In this section, we explore two algorithms to approximate nonsymmetric matrices under the basic framework $A \approx \hat{A}_{r} = AX(Y^{T}AX)^{\dagger}(A^{T}Y)^{T}$. There is no doubt that the GN-c is more attractive in comparison with GN-r\&c and GN. More comparisons of their performances are illuminated in the following numerical experiments, in which different singular value decaying speeds of $A$ will be taken. For GN-r\&c and GN-c, we perform the \textit{economy-size} QR factorization frequently, which operates on the whole matrix at once with an expensive cost. For improving the performances of the methods, we can choose the column-pivoted QR factorization by applying groups of Householder reflectors \cite{Golub} or other randomized rank revealing QR factorizations in \cite{demmel2015,gu2020}.

    \section{Error analysis}
   In \cite{gittens2013}, for a symmetric positive semidefinite matrix $A \succeq 0$, its rank-$r$ approximation is obtained by the form $AQ(Q^{T}AQ)^{\dagger}(AQ)^{T}$, where $Q = \mathrm{orth}(A\Omega)$ and $\Omega$ is a sketching matrix. The resulting accuracy is considerably better than that of rSVD.
   For general matrices, the GN in \cite{Naka2020} shows its numerical stability and fast speed. But in accuracy, the GN does not outperform rSVD and the extended methods, including GN-r\&c and GN-c with better numerical performances as we predicted. Meanwhile, the GN-r\&c is not recommend as a practical method for its expensive computational cost. In this section, we will analyse the approximation accuracy with the form $\lVert A - \hat{A}_{r}\lVert_{F}$ and focus on the comparisons among GN, rSVD, and GN-c.
   \subsection{Projections and key properties}
    At the beginning, it is necessary to note some projections onto different spaces in these methods mentioned above. In order to facilitate the comparisons and illustrations, we denote related projectors by one single unified form as follows
   \begin{equation}
    \mathcal{P}_{X,Y} \equiv X(Y^{T}X)^{\dagger}Y^{T}.
   \end{equation}

   As we all know, an orthogonal projector is a Hermitian matrix $\mathcal{P}$, which is completely determined by its range. For a matrix $X$, we can write the orthogonal projection $\mathcal{P}_{X}$ with range($\mathcal{P}_{X}$) = range($X$). And the orthogonal projection can be expressed explicitly as
   \begin{equation}
    \mathcal{P}_{X} \equiv X(X^{T}X)^{\dagger}X^{T},
   \end{equation}
   when $X = Y$ in (5), then $\mathcal{P}_{X,Y} = \mathcal{P}_{X,X} = \mathcal{P}_{X}$ \cite{Naka2020}. The expressions of (5) and (6) will be used in the following analysis about the projectors of rSVD, GN, and GN-c, which we refer to as $\mathcal{P}_{rSVD}, \mathcal{P}_{GN},$ and $\mathcal{P}_{GN-c}$, respectively.

   When considering the rSVD, we form the approximation $\hat{A}_{r} = QQ^{T}A$, where $Q = \mathrm{orth}(A\Omega)$, so that range($A\Omega$) = range(Q). Because the range($A\Omega$) is consistent with the range($\mathcal{P}_{A\Omega}$), we can express the orthogonal projector of rSVD as
   \begin{equation}
    \mathcal{P}_{rSVD} \equiv \mathcal{P}_{A\Omega} = A\Omega((A\Omega)^{T}A\Omega)^{\dagger}(A\Omega)^{T}.
   \end{equation}
   And the rank-$r$ approximation via rSVD can be expressed as
   \begin{equation}
    \hat{A}_{rSVD} = \mathcal{P}_{rSVD}A = \mathcal{P}_{A\Omega}A.
   \end{equation}

   Meanwhile, the approximation form of GN is $\hat{A}_{GN} = AX(Y^{T}AX)^{\dagger}Y^{T}A$, and we obtain two key oblique projections
   \begin{equation}
    \mathcal{P}_{AX,Y} \equiv AX(Y^{T}AX)^{\dagger}Y^{T},\quad
    \mathcal{P}_{X,A^{T}Y} \equiv X(Y^{T}AX)^{\dagger}Y^{T}A,
   \end{equation}
   which are projected onto the column space of $AX$ and row space of $Y^{T}A$, respectively. As explained in \cite{Naka2020}, we have\\
   \begin{equation}
   \hat{A}_{GN} = \mathcal{P}_{AX,Y}A = A\mathcal{P}_{X,A^{T}Y} = \mathcal{P}_{AX,Y}A\mathcal{P}_{X,A^{T}Y},
   \end{equation}
   which is easy to be verified and the formula above is also applicable to the orthogonal projection of rSVD.

   As for GN-c method, the rank-$r$ approximation of $A$ is $\hat{A}_{GN\mbox{-}c} = A\hat{Q}(Q_{1}^{T}A\hat{Q})^{\dagger}(A^{T}Q_{1})^{T}$, where $Q_{1} = \mathrm{orth}(AX)$ and $\hat{Q} = \mathrm{orth}(A^{T}Q_{1})$, which merely generates one Gaussian matrix $X$. The projector of GN-c method is
   \begin{equation}
    \mathcal{P}_{GN\mbox{-}c} \equiv \mathcal{P}_{\hat{Q},A^{T}Q_{1}} = \hat{Q}(Q_{1}^{T}A\hat{Q})^{\dagger}(A^{T}Q_{1})^{T}.
   \end{equation}

   In the following, we will focus on the rSVD, GN, and GN-c methods and investigate their theoretical error bounds. The analysis will be based on the original approximation form $AX(Y^{T}AX)^{\dagger}(A^{T}Y)^{T}$ and these significant projections mentioned above. Before we do, it is necessary to summarize some key properties of projections as follows.
   \begin{itemize}
    \item[(\romannumeral 1)] $\lVert I - \mathcal{P}\lVert_{2} = \lVert \mathcal{P}\lVert_{2}$, which holds for any projection $\mathcal{P}$ s.t. $\mathcal{P} = \mathcal{P}^{2}$ \cite{szyld2006};
    \item[(\romannumeral 2)] For oblique projections of $\mathcal{P}_{AX,Y}$ and $\mathcal{P}_{X,A^{T}Y}$ with $Y^{T}AX$ having full column rank, $(I - \mathcal{P}_{AX,Y})AX = 0$ and $Y^{T}A(I - \mathcal{P}_{X,A^{T}Y}) = 0$;
    \item[(\romannumeral 3)] For the orthogonal projection $\mathcal{P}$, we have $\mathcal{P}^{T} = \mathcal{P}$ and $\lVert I - \mathcal{P}\lVert_{2} = \lVert \mathcal{P}\lVert_{2} = 1$.
   \end{itemize}
   \subsection{Approximation accuracy}
   In this subsection, firstly, some useful conclusions of the random matrix theory have to be presented as follows.
   \begin{proposition}\textnormal{(Gaussian matrix products)} \cite{HMT2011}
    Fix a matrix $S\in\mathbb{R}^{(m-k)\times(m-k)}$ and consider independent Gaussian matrices $G\in\mathbb{R}^{(m-k)\times r}$ and $H\in\mathbb{R}^{k\times r}$, where $r\geq k.$ Then, for any $u, t \geq 0$,\\
    $$\mathbb{P}\left\lbrace \lVert SGH^{\dagger}\lVert_{F}^{2}>\frac{utk}{r-k+1}\lVert S\lVert_{F}^{2} \right\rbrace\leq e^{-(u-2)/4}+\sqrt{\pi k}(t/e)^{-(r-k+1)/2}. $$\\
    Additionally, if $r\geq k+2,$\\
    $$\mathbb{E}\lVert SGH^{\dagger}\lVert_{F}^{2} = \frac{k}{r-k-1}\lVert S\lVert_{F}^{2}.$$
   \end{proposition}

   Another useful fact of rectangular random matrix should be noted as well. For random matrices with entries from the standard normal distribution $\mathcal{N}(0, 1)$, (i.e. the Gaussian matrices), the calssical Marchenko-Pastur (M-P) rule \cite{MP} shows that a random matrix $\Omega\in\mathbb{R}^{m\times n} (m\geq n)$ has singular values supported in the interval $[\sqrt{m} - \sqrt{n}, \sqrt{m} + \sqrt{n}]$, which identifies the precise limiting distribution.
    \subsubsection{Analysis of rSVD \& GN}
    Some key steps of following numerical analysis are inspired by the analysis in \cite{sorensen2016}. It is necessary to note that rectangular random matrices are well-conditioned, which will be used repeatedly in the forthcoming section.

    As we can see, when choosing $Y = AX$, the projection $\mathcal{P}_{AX,Y}$ in (9) becomes an orthogonal projection,
    $$\mathcal{P}_{AX,Y} = AX((AX)^{T}AX)^{\dagger}(AX)^{T},$$
    which is consistent with the projection (7) of rSVD. So the error matrix of rSVD can be unified into the error form of GN, which can be expressed as
    $$E = A - AX(Y^{T}AX)^{\dagger}Y^{T}A = (I - \mathcal{P}_{AX,Y})A = (I - \mathcal{P}_{AX,Y})A(I - \mathcal{P}_{X,A^{T}Y}),$$
    where the third equality can be proved using the formula (10).

    Taking the property (\romannumeral 2) of the oblique projection into consideration, we can derive another equivalent form
    \begin{equation}
        E = (I - \mathcal{P}_{AX,Y})A = (I - \mathcal{P}_{AX,Y})A(I - X\Pi_{X}),
    \end{equation}
    for any matrix $\Pi_{X} \in \mathbb{R}^{r \times n}$. The second equality holds because $(I - \mathcal{P}_{AX,Y})A(I - X\Pi_{X}) = (I - \mathcal{P}_{AX,Y})A - (I - \mathcal{P}_{AX,Y})AX\Pi_{X}$, where $(I - \mathcal{P}_{AX,Y})AX\Pi_{X} = 0$ with $Y^{T}AX$ having full column rank.

    Nextly, let us focus on the error bound of rSVD, which will be used to compare with the error bounds of GN and GN-c. And we define the error of rSVD as $E_{rSVD} = A - \hat{A}_{rSVD}$. The following procedures are useful for the proofs of Theorem 4.2 and Theorem 4.3. According to (12), we have
    \begin{equation}
        E_{rSVD} = A - \hat{A}_{rSVD} = (I - \mathcal{P}_{A\Omega})A = (I - \mathcal{P}_{A\Omega})A(I - X\Pi_{X}),
    \end{equation}
    and let $\Pi_{X}$ take the form $(W^{T}X)^{\dagger}W^{T}$, where $W \in \mathbb{R}^{n \times k}$ is arbitrary and $k \leq r$ is a target rank.

    For the orthogonal projector $\mathcal{P}_{A\Omega}$, $\lVert \mathcal{P}_{A\Omega}\lVert_{2} = \lVert I - \mathcal{P}_{A\Omega}\lVert_{2} = 1$. Then we obtain
    \begin{equation}
        \lVert E_{rSVD}\lVert_{2} = \lVert (I - \mathcal{P}_{A\Omega})A\lVert_{2} = \lVert (I - \mathcal{P}_{A\Omega})A(I - X\Pi_{X})\lVert_{2} \leq \lVert A(I - X\Pi_{X})\lVert_{2}.
    \end{equation}

    Now let $W$ consist of the first $k$ leading left singular vectors of $A^{T}$, and it is obvious that $X\Pi_{X} = X(W^{T}X)^{\dagger}W^{T}$ is an oblique projection onto the subspace of $X$. Assuming that $W^{T}X$ has full row-rank, thus $W^{T}X(W^{T}X)^{\dagger} = I$, so that $W^{T}(I - X\Pi_{X}) = 0$.

    Furthermore, noting that $W_{\perp}W_{\perp}^{T} = I - WW^{T}$, thus $AW_{\perp}W_{\perp}^{T} = \Sigma$, where $\Sigma = \mathrm{diag}(\sigma_{k+1}, \cdots , \sigma_{n}).$ Then we have
    $$\lVert A(I - X\Pi_{X})\lVert_{2} = \lVert A(I - WW^{T})(I - X\Pi_{X})\lVert_{2} = \lVert AW_{\perp}W_{\perp}^{T}(I - X\Pi_{X})\lVert_{2}.$$

     By the equality $\lVert AW_{\perp}W_{\perp}^{T}\lVert_{2} = \lVert \Sigma\lVert_{2}$, we can obtain
    \begin{equation}
        \mathbb{E}\lVert E_{rSVD}\lVert_{2} \leq \lVert A(I - X\Pi_{X})\lVert_{2} \leq \lVert\Sigma\lVert_{2}\lVert I - X\Pi_{X}\lVert_{2} = \lVert\Sigma\lVert_{2}\lVert X\Pi_{X}\lVert_{2}.
    \end{equation}

     It is easy to verify that $\lVert X\Pi_{X}\lVert_{2} = \lVert X(W^{T}X)^{\dagger}W^{T}\lVert_{2} = \lVert X(W^{T}X)^{\dagger}\lVert_{2}$. Thus by the M-P rule \cite{MP}, supposing that $X\in\mathbb{R}^{n\times r}$ is Gaussian and so is $W^{T}X\in\mathbb{R}^{k\times r}$ $(k < r)$, we have $\mathbb{E}\lVert (W^{T}X)^{\dagger}\lVert_{2} \approx \frac{1}{(\sqrt{r} - \sqrt{k})}$ and $\mathbb{E}\lVert X\lVert_{2} \approx \sqrt{m} + \sqrt{r}$ with high probability. Thus
     \begin{equation}
        \mathbb{E}\lVert E_{rSVD}\lVert_{2} \lesssim \dfrac{\sqrt{m} + \sqrt{r}}{\sqrt{r} - \sqrt{k}}\lVert \Sigma\lVert_{2}.
     \end{equation}

     We can find that rSVD is optimal to within the oversampling $r - k$ and the factor $(\sqrt{m} + \sqrt{r})/(\sqrt{r} - \sqrt{k})$. The analysis above is merely a more detailed proof of the exisiting conclusion in \cite{Naka2020}. Either the error bound in the spectral norm or the following one in the Frobenius norm is useful to estimate the error of GN-c.

    And moving forward, we focus on the error bound in the Frobenius norm of rSVD and GN.
    \begin{theorem}\cite{HMT2011}
        Fix a general matrix $A\in\mathbb{R}^{m\times n}$ and a target rank $k \geq 1$, and let $A_{k}$ be an optimal rank-$k$ approximation that arises from a $k$-truncated singular value decomposition of $A$. Then rSVD with $r \geq k + 2$ Gaussian initialization vectors generates a rank-$r$ approximation $\hat{A}_{rSVD}$ that satisfies
        \begin{equation}
                \mathbb{E}\lVert E_{rSVD}\lVert_{F} = \mathbb{E}\lVert A - \hat{A}_{rSVD}\lVert_{F} \leq \sqrt{\mathbb{E}\lVert A - \hat{A}_{rSVD} \lVert_{F}^{2}} = \sqrt{1 + \dfrac{k}{r - k -1}}\lVert A - A_{k}\lVert_{F}.
        \end{equation}
    \end{theorem}

    Turn to the GN method and define the error matrix of GN as $E_{GN} = A - \hat{A}_{GN}$. We again start with (12), and a difference from rSVD is that $\mathcal{P}_{AX,Y}$ is an oblique projector. Hence we will show the error bound in the Frobenius norm directly in Theorem 4.3.
    \begin{theorem}\cite{Naka2020}
        Fix the general matrix $A\in\mathbb{R}^{m\times n}$, a target rank $k \geq 1$ and the oversampling parameter $l$. Let $A_{k}$ be an optimal rank-$k$ approximation of $A$. Then a rank-$r$ approximation via GN $\hat{A}_{GN} \in \mathbb{R}^{m\times n}$ satisfies
        \begin{equation}
            \mathbb{E}\lVert E_{GN}\lVert_{F} = \mathbb{E}\lVert A - \hat{A}_{GN}\lVert_{F} \leq \sqrt{\left( 1+ \frac{r + l}{l - 1}\right)\left( 1 + \dfrac{k}{r - k -1}\right) }\lVert A - A_{k}\lVert_{F}.
        \end{equation}
    \end{theorem}

    The error bounds in Theorem 4.2 and Theorem 4.3 show that rSVD is optimal to within a small oversampling, while the GN is susceptible to the oversampling $l$, especially when the singular value decaying speed is not fast. More details of the related proofs for Theorem 4.2 and Theorem 4.3 can be found in \cite{HMT2011,Naka2020,Tropp2023}.
    \subsubsection{Analysis of GN-c}
    Nextly, we will utilize the insights of randomized subspace iteration \cite{HMT2011,gu2015} and the related projections to analyse the error bound of GN-c method. Meanwhile the error of GN-c is defined as $E_{GN\mbox{-}c} = A - \hat{A}_{GN\mbox{-}c}$.

    \begin{proposition}\cite{HMT2011}
        Let $\mathcal{P}$ be an orthogonal projector, and let $V$ be an arbitrary matrix with compatible dimensions. For each positive number $q$,
        \begin{equation}
            \lVert \mathcal{P}V\lVert \leq \lVert \mathcal{P}(VV^{T})^{q}V\lVert^{1/(2q+1)}.
        \end{equation}
    \end{proposition}

    To begin, some notations must be introduced, and the key error bound will be established afterward. Let $A$ be an $m\times n$ matrix that has an SVD $A = U\Sigma V^{T}$. To simplify the analysis, it is appropriate to partiton the SVD as follows:
    \begin{equation}
        A = U
        \begin{pmatrix}
            \Sigma_{1} &  \\
                          & \Sigma_{2}\\
        \end{pmatrix}
    \begin{pmatrix}
        V_{1}^{T}\\
        V_{2}^{T}\\
    \end{pmatrix}.
    \end{equation}
     We can see that the left unitary matrix $U$ does not play a significant role in the analysis.
     Then let $\Omega$ be an $n\times r$ Gaussian matrix, where we assume $r \geq k+2$, and partion $V^{T}\Omega$ into
     \begin{equation}
        \Omega_{1} = V_{1}^{T}\Omega,\quad
        \Omega_{2} = V_{2}^{T}\Omega.
     \end{equation}

    With the notations above, we continue the following analysis. The projections of GN-c and the randomized subspace iteration methods in \cite{bjark,HMT2011} and the connections among them play a significant role in the analysis. Related conclusions below in Theorem 4.6 and Lemma 4.7 are helpful to prove the error bound of GN-c.
    \begin{lemma} \cite{HMT2011}.
        Let $A$ be an $m\times n$ matrix with SVD $A = U\Sigma V^{T}$. Choose a Gaussian matrix $\Omega\in\mathbb{R}^{n\times r}$ and construct the sample matrix $A\Omega$. Partion $\Sigma$ as specified in \textnormal{(20)}, and define $\Omega_{1}$ and $\Omega_{2}$ as \textnormal{(21)}. Assuming that $\Omega_{1}$ has full row rank, the approximation error satisfies
        \begin{equation}
         \lVert (I - \mathcal{P}_{A\Omega})A\lVert^{2} \leq \lVert \Sigma_{2}\lVert^{2} + \lVert \Sigma_{2}\Omega_{2}\Omega_{1}^{\dagger}\lVert^{2},
        \end{equation}
    where $\lVert  \cdot\lVert$ denotes either the spectral norm or the Frobenius norm.
    \end{lemma}

    Applying Proposition 4.1 directly, we have $\mathbb{E}\lVert \Sigma_{2}\Omega_{2}\Omega_{1}^{\dagger}\lVert^{2}_{F} = \dfrac{k}{r-k-1}\lVert \Sigma_{2}\lVert^{2}_{F}$, so that
    $$\mathbb{E}\lVert (I - \mathcal{P}_{A\Omega})A\lVert_{F} \leq \sqrt{1 + \dfrac{k}{r - k -1}}\lVert \Sigma_{2}\lVert_{F},$$
    which is useful to prove Theorem 4.6.

    And then when considering to approximate $A$ via the standard randomized subspace iteration \cite{HMT2011}, we fix a positive integer $q$, and set
   \begin{equation}
    B = (AA^{T})^{q}A = U\Sigma^{2q+1}V^{T}.
   \end{equation}
    Then we generate the basis $Q_{c}$ to approximate the range of $B\Omega$ where $Q_{c} = \mathrm{orth}((AA^{T})^{q}A\Omega)$, and the rank-$r$ approximation of $A$ via standard randomized subspace iteration can be expressed as
        \begin{equation}
            \hat{A}_{r} = \mathcal{P}_{B\Omega}A = \mathcal{P}_{(AA^{T})^{q}A\Omega}A = Q_{c}Q_{c}^{T}A.
        \end{equation}

     The following result in Theorem 4.6 describes the error bound $\lVert (I-\mathcal{P}_{B\Omega})A\lVert$ in the Frobenius norm.
    \begin{theorem}
        Let $A$ be an $m\times n$ matrix with SVD $A = U\Sigma V^{T}$, and choose a Gaussian matrix $\Omega\in\mathbb{R}^{n\times r}$. Fix an integer $q\geq 0$, and form $B\Omega = (AA^{T})^{q}A\Omega$. Partion $\Sigma$ as specified in \textnormal{(20)}, and define $\Omega_{1}$ and $\Omega_{2}$ as \textnormal{(21)}. Assuming that $\Omega_{1}$ has full row rank, then
        \begin{equation}
            \mathbb{E}\lVert (I-\mathcal{P}_{B\Omega})A\lVert_{F} \leq \left( 1+\dfrac{k}{r-k-1}\right) ^{1/(4q+2)}\lVert \Sigma_{2}^{2q+1}\lVert_{F}^{1/(2q+1)}.
        \end{equation}
    \end{theorem}
    \textit{Proof.} As a direct consequence of Proposition 4.4, we conclude that
    $$\lVert (I-\mathcal{P}_{B\Omega})A\lVert \leq \lVert (I-\mathcal{P}_{B\Omega})B\lVert_{F}^{1/(2q+1)} = \lVert (I-\mathcal{P}_{(AA^{T})^{q}A\Omega})(AA^{T})^{q}A\lVert^{1/(2q+1)},$$
    then applying Lemma 4.5 and (23), we have
   \begin{equation}
     \lVert (I-\mathcal{P}_{(AA^{T})^{q}A\Omega})(AA^{T})^{q}A\lVert^{2/(2q+1)} \leq (\lVert \Sigma_{2}^{2q+1}\lVert^{2} + \lVert \Sigma_{2}^{2q+1}\Omega_{2}\Omega_{1}^{\dagger} \lVert^{2})^{1/(2q+1)}.
   \end{equation}
   Using the Frobenius norm and applying Proposition 4.1, we have
   $$(\lVert \Sigma_{2}^{2q+1}\lVert^{2}_{F} + \mathbb{E}\lVert \Sigma_{2}^{2q+1}\Omega_{2}\Omega_{1}^{\dagger} \lVert^{2}_{F})^{1/(2q+1)} \leq \left[ \left( 1+\dfrac{k}{r-k-1}\right) \lVert \Sigma_{2}^{2q+1}\lVert_{F}^{2}\right] ^{1/(2q+1)}.$$
   It follows that
   $$\mathbb{E}\lVert (I-\mathcal{P}_{B\Omega})A\lVert_{F} \leq \mathbb{E}\lVert (I-\mathcal{P}_{B\Omega})B\lVert_{F}^{1/(2q+1)} \leq \left( 1+\dfrac{k}{r-k-1}\right) ^{1/(4q+2)}\lVert \Sigma_{2}^{2q+1}\lVert_{F}^{1/(2q+1)}. $$

    \qquad \qquad \qquad \qquad \qquad \qquad \qquad \qquad \qquad \qquad \qquad  \qquad \qquad \qquad  \qquad\qquad \qquad$\square$

    In \cite{bjark}, unlike the standard form of randomized subspace iteration in (24), the approximation $\hat{A}_{r} = \mathcal{P}_{(AA^{T})^{q}A\Omega}A$ is replaced with the closely related form
        $$\hat{A}_{r} = A\mathcal{P}_{(A^{T}A)^{q}\Omega}.$$

    The basic idea of the generalized subspace iteration method presented above is that the standard subspace iteration process can be halted halfway through an iteration to give an orthonormal matrix $Q_{r}$ which approximates the row space of $A$ where $Q_{r} = \mathrm{orth}((A^{T}A)^{q}\Omega)$ \cite{bjark}. This corresponds to evaluating an orthonormal basis to approximate the range of $(A^{T}A)^{q}\Omega$. Then we have
    \begin{equation}
        \hat{A}_{r} = A\mathcal{P}_{(A^{T}A)^{q}\Omega} = AQ_{r}Q_{r}^{T},
    \end{equation}
     In (24), the standard subspace iteration method generates $Q_c$ to approximate the column space of $A$ performing $2q+1$ power iterations. Similarly, the form in (27) of generalized subspace iteration method depends on $2q$ power iterations to form $Q_{r}$, which approximates the row space of $A$.

    Inspired by the relation between the standard and generalized subspace iteration methods, it is convenient to apply the result of Theorem 4.6 to derive the following inequality.
    \begin{equation}
        \mathbb{E}\lVert A(I - \mathcal{P}_{(A^{T}A)^{q}\Omega})\lVert_{F} \leq \left( 1+\dfrac{k}{r-k-1}\right) ^{1/4q}\lVert \Sigma_{2}^{2q}\lVert_{F}^{1/2q}.
    \end{equation}

    The error bound in (28) is what we can expect based on Theorem 4.6. That is, the error decays exponentially with the iteration number $q$. Nextly, let us make full use of the conclusions above and focus on the error bound of GN-c. Also define the error matrix of GN-c as $E_{GN\mbox{-}c} = A - \hat{A}_{GN\mbox{-}c}.$
    \begin{lemma}
        Let $\mathcal{P}_{X,Y} = X(Y^{T}X)^{\dagger}Y^{T}$ be a projector. If range(X) = range(Y), the projection $\mathcal{P}_{X,Y}$ is orthogonal.
    \end{lemma}

    The rank-$r$ approximation of $A$ via GN-c can be equally written as $\hat{A}_{GN\mbox{-}c} = A\hat{Q}(Q_{1}^{T}A\hat{Q})^{\dagger}(A^{T}Q_{1})^{T}$, where $Q_{1} = \mathrm{orth}(AX)$ and $\hat{Q} = \mathrm{orth}(A^{T}Q_{1})$.
    It is obvious that the range of $\hat{Q}$ is consistent with the range of $A^{T}Q_{1}$, so that according to Lemma 4.7, the projector (11) of GN-c method
    $$\mathcal{P}_{GN\mbox{-}c} = \mathcal{P}_{\hat{Q},A^{T}Q_{1}} = \hat{Q}(Q_{1}^{T}A\hat{Q})^{\dagger}(A^{T}Q_{1})^{T} = \hat{Q}(\hat{Q}^{T}\hat{Q})^{\dagger}\hat{Q}^{T}$$
    is an orthogonal projection.

    Owning to $\hat{Q} = \mathrm{orth}(A^{T}AX)$, the rank-$r$ approximation by GN-c can be expressed equivalently as
    $$\hat{A}_{GN\mbox{-}c} = A\mathcal{P}_{GN\mbox{-}c} =  A\mathcal{P}_{\hat{Q}} = A\mathcal{P}_{A^{T}AX}.$$

    So for GN-c method, applying the result of (28) with $q = 1$, we obtain the error bound of GN-c method
    $$\mathbb{E}\lVert E_{GN\mbox{-}c}\lVert_{F} = \mathbb{E}\lVert A - \hat{A}_{GN\mbox{-}c}\lVert_{F} \leq \left( 1+\dfrac{k}{r-k-1}\right) ^{1/4}\lVert \Sigma_{2}^{2}\lVert_{F}^{1/2}.$$

    Additionally, for the spectral norm, we have
    $$\mathbb{E}\lVert E_{GN\mbox{-}c}\lVert_{2} = \mathbb{E}\lVert A - \hat{A}_{GN\mbox{-}c}\lVert_{2} \leq \left[\left( 1 + \sqrt{\dfrac{k}{r - k -1}}\right) \sigma_{k + 1}^{2} + \dfrac{e\sqrt{r}}{r-k}\left( \sum\limits_{i > k}\sigma_{i}^4\right) ^{1/2} \right]^{1/2},$$
    which can be proved by the techniques in \cite[Thm. 4.1]{bjark}.

    As mentioned above, the related analysis and comparisons of the error bounds for GN, rSVD and GN-c have been introduced. The procedures of obtaining the error bound of GN-c have been shown in detail. Then we summarize the error bounds of GN-c method in Theorem 4.8.

    \begin{theorem}
        Fix a matrix $A \in \mathbb{R}^{m\times n}$ with singular values $\sigma_{1} \geq \sigma_{2} \geq \sigma_{3} \geq \cdots.$ And its rank-$r$ approximation is $\hat{A}_{GN\mbox{-}c}$ via GN-c.
        Then for any $r \geq k + 2$, where the target rank $k \geq 2$, the error bound in the spectral norm satisfies
        $$\mathbb{E}\lVert E_{GN\mbox{-}c}\lVert_{2} \leq \left[\left( 1 + \sqrt{\dfrac{k}{r - k -1}}\right) \sigma_{k + 1}^{2} + \dfrac{e\sqrt{r}}{r-k}\left( \sum\limits_{i > k}\sigma_{i}^4\right) ^{1/2} \right]^{1/2},$$
        additionally in the Frobenius norm,
        $$\mathbb{E}\lVert E_{GN\mbox{-}c}\lVert_{F} \leq \left( 1+\dfrac{k}{r-k-1}\right) ^{1/4}\lVert \Sigma_{2}^{2}\lVert_{F}^{1/2}.$$
    \end{theorem}

     Meanwhile, the results in numerical experiments are helpful to intuitively sense the differences between all the methods above and verify the error bounds we have derived. GN-c method does work well in practice as we expect.
    \section{Numerics}\setcounter{footnote}{0}\renewcommand{\thefootnote}{\arabic{footnote}}
    In section 3, we have briefly displayed the performances of rSVD, GN, and GN-r\&c in terms of relative errors and computational speeds. It is obvious that the singular value decaying speed of $A$ plays a significant role in the performances of the above algorithms for low-rank matrix approximation. In this section, we will test different data examples to explore the performances of these methods and verify the analysis results. Also all numerical experiments have been performed in MATLAB (version 2022a) on a MacBook Pro with a 2.3 GHz Intel Core i7 processor with four cores.

    The experiments consist of the synthetic examples and dataset examples\footnote[1]{https://math.nist.gov/MatrixMarket/}\footnote[2]{https://sparse.tamu.edu} and all the matrices are nonsymmetric with different singular spectrum decaying speeds. The results of numerical experiments suggest that GN-c outperforms other approaches in practice and confirm the error bounds in Theorem 4.8 meanwhile.

     We will measure the relative errors with the form
     $$\dfrac{\lVert A - \hat{A}_{r}\lVert_{F}}{ \lVert A\lVert_{F}},$$
     to examine the accuracy of these algorithms. And the runtimes of performing the program for the methods, including GN, rSVD, GN-r\&c, and GN-c are recorded to show their computational speeds.
    \subsection{Synthetic examples}
    For demonstrating the efficiency of our proposed methods adequately, we utilize four synthetic matrices, whose singular value decaying speeds are different.
    We construct the synthetic test matrices as follows.
    \begin{itemize}
        \item[(1)] The synthetic matrices take the same form $A = M\Lambda N^{T}$, where $M\in\mathbb{R}^{m\times m}$ and $N\in\mathbb{R}^{m\times m}$ are the Q-factors of two different independent Gaussian matrices, respectively, and $\Lambda$ is diagonal and its singular spectrum decays with four different decaying speeds.
        \item[(2)] The singular value decaying speeds of $\Lambda$ is taken as follows
        \begin{itemize}
            \item[(a)] Polynomially decaying spectrum: $\Lambda = \mathrm{diag}(1, 2^{-p}, 3^{-p}, \cdots , m^{-p}) \in \mathbb{R}^{m\times m}$
            \begin{itemize}
                \item Slow polynomial decay (PolyDecaySlow): $p = 1.$
                \item Fast polynomial decay (PolyDecayFast): $p = 2.$
            \end{itemize}
            \item[(b)] Exponentially decaying spectrum: $\Lambda = \mathrm{diag}(1, 10^{-q}, 10^{-2q}, \cdots , 10^{-(m-1)q}) \in \mathbb{R}^{m\times m}$
            \begin{itemize}
                \item Slow exponential decay (ExpDecaySlow): $q = 0.125.$
                \item Fast exponential decay (ExpDecayFast): $q = 0.25.$
            \end{itemize}
        \end{itemize}
    \end{itemize}
    Here, the values of $p$ and $q$ are chosen to control the rates of different singular value decaying speeds.

    Note that in the plots, the Algorithm 2-1 \cite{Naka2020} and Algorithm 2-2 \cite{HMT2011} are referred to as GN and rSVD, and the algorithms proposed by us in section 3 are called GN-r\&c and GN-c. In every experiment, we compare the relative errors and speeds with varying rank-$r$ to demonstrate the anvantages and disadvantages of GN, rSVD, GN-r\&c, and GN-c. Then we start to illustrate the performances of all the methods mentioned up to now via the numerical experiments.
    \subsubsection{Synthetic examples with polynomially decaying spectrum}
\begin{figure} [H]
        \subfigure[\textit{Synthetic matrix (PolyDecaySlow)}]{
            \includegraphics[width=8cm,height=6cm]{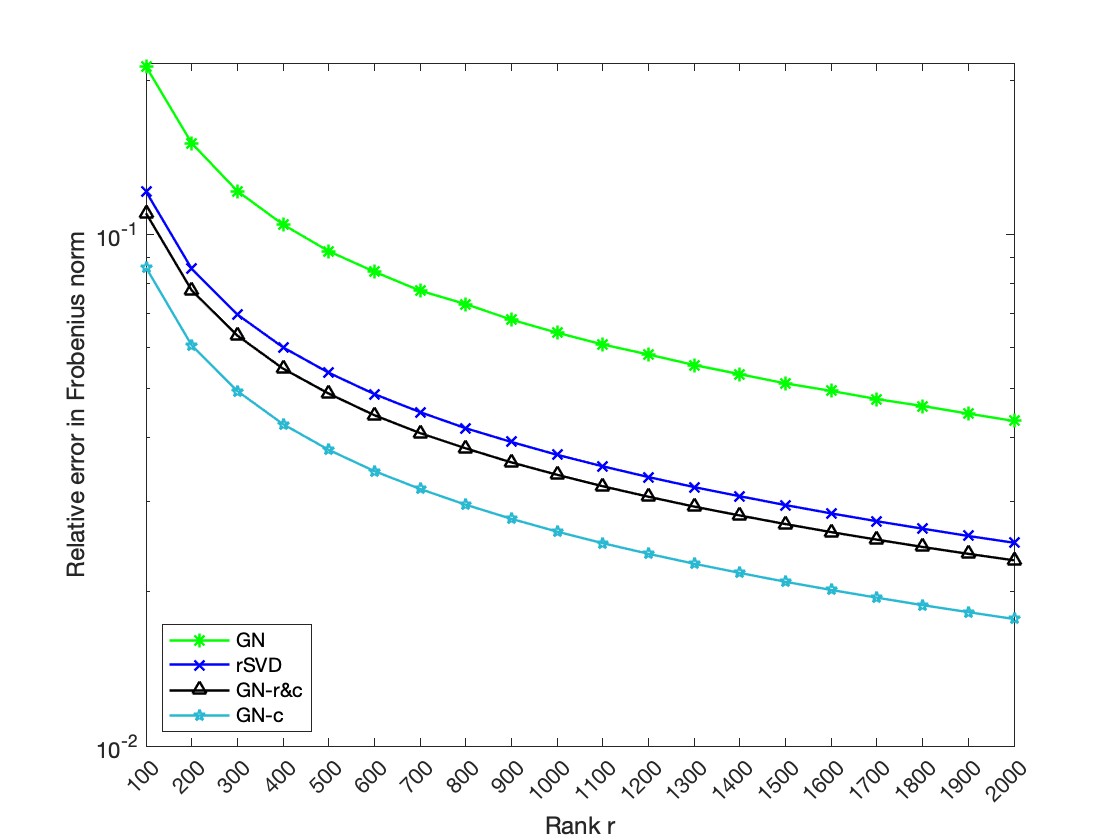} \label{Fig.sub.sp}
        }
        \hspace{2mm}
        \subfigure[\textit{Synthetic matrix (PolyDecaySlow)}]{
            \includegraphics[width=8cm,height=6cm]{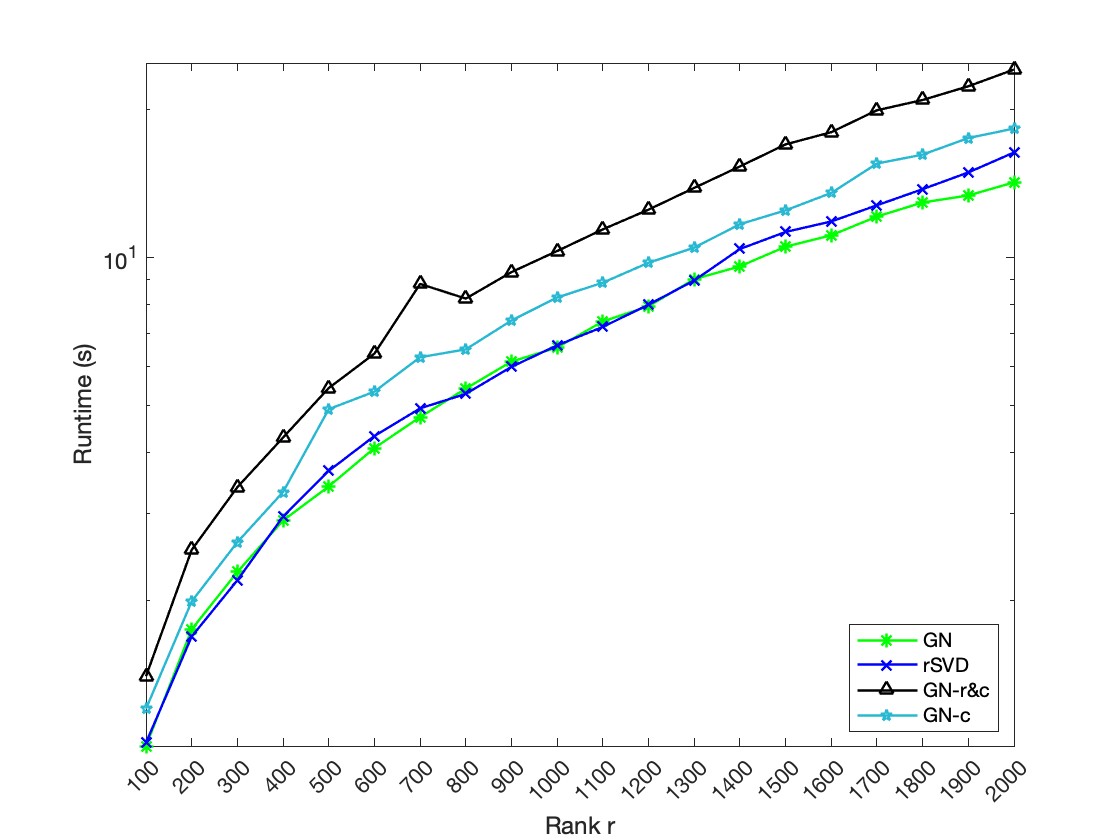} \label{Fig.sub.spt}
        }

        \subfigure[\textit{Synthetic matrix (PolyDecayFast)}]{
            \includegraphics[width=8cm,height=6cm]{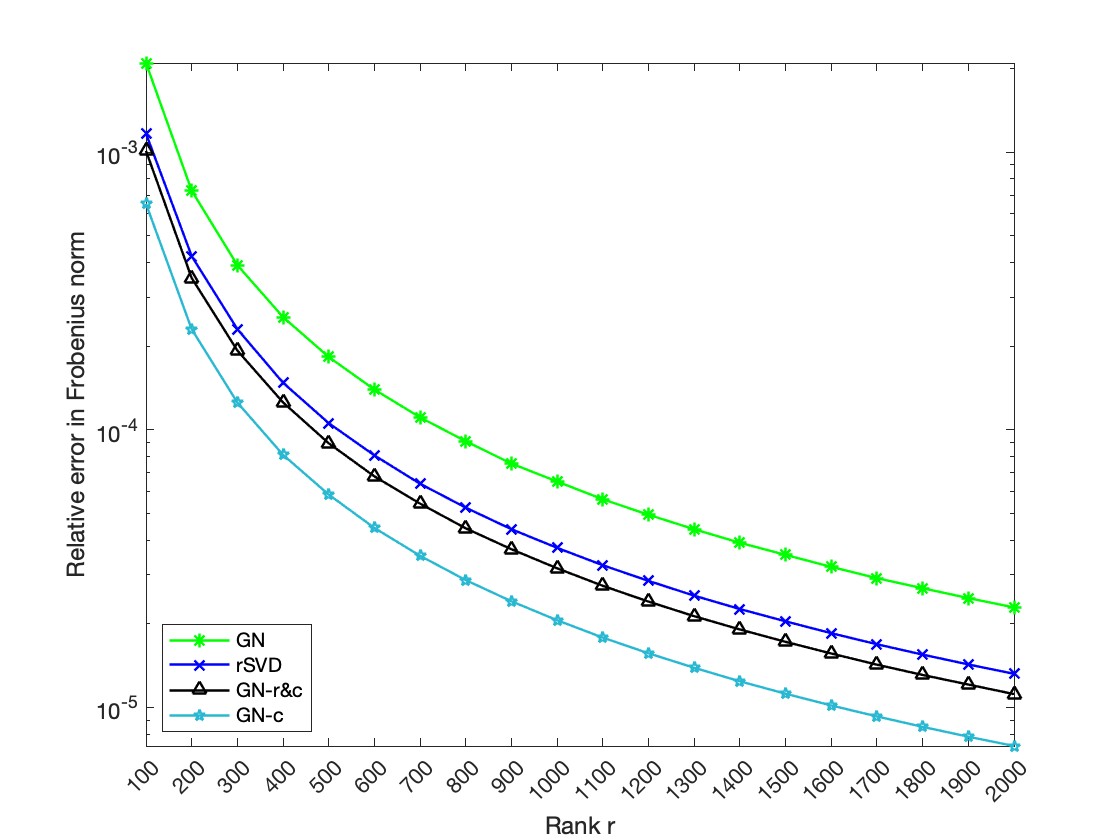} \label{Fig.sub.fp}
        }
        \hspace{2mm}
        \subfigure[\textit{Synthetic matrix (PolyDecayFast)}]{
            \includegraphics[width=8cm,height=6cm]{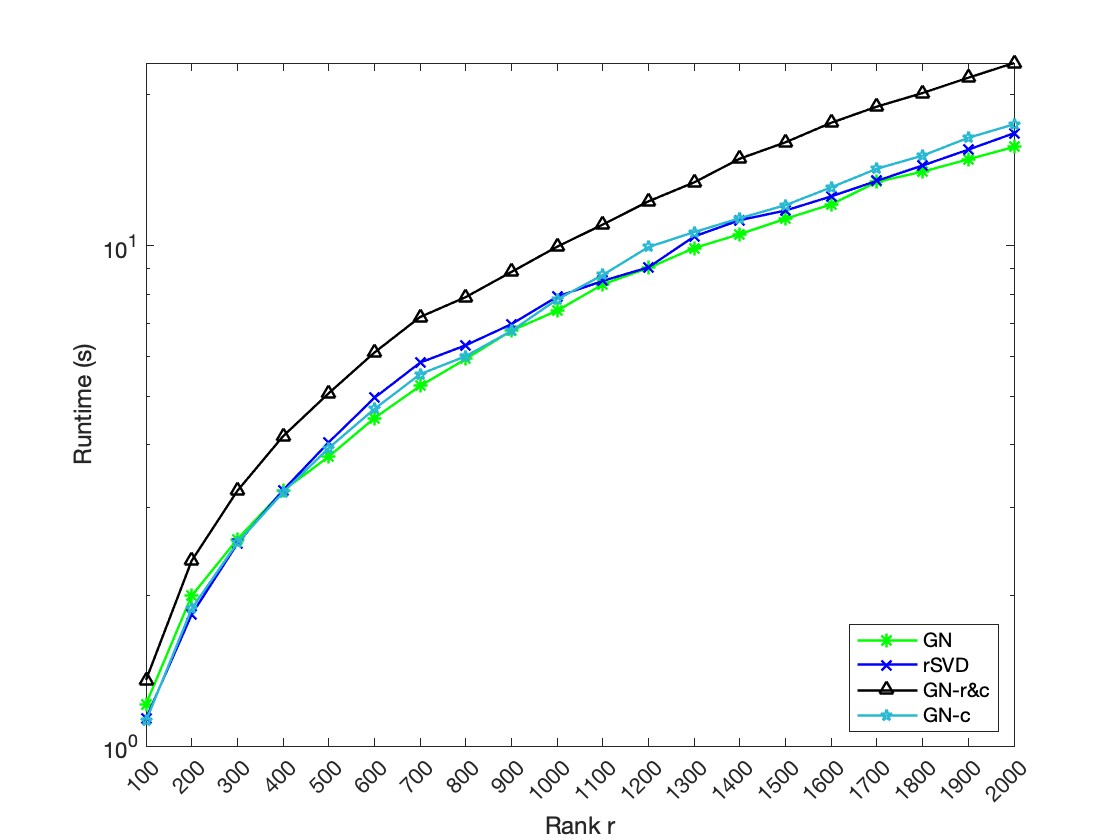} \label{Fig.sub.fpt}
        }
        \caption{\textit{\textbf{Relative error and runtime comparisons.}  The relative error and runtime of approximating the matrix with slow polynomial decaying speed are shown in (\romannumeral 1) and (\romannumeral 2). And the relative error and runtime of approximating the matrix with fast polynomial decaying speed are presented in (\romannumeral 3) and (\romannumeral 4).}}\label{Fig.main.SP}
     \end {figure}

    In Fig.\ref{Fig.main.SP}, we fix the dimension parameters of test matrices $m = n = 10000$.
    The singular value decaying speeds of the synthetic matrices are slow polynomial ($p = 1$) and fast polynomial ($p = 2$), respectively. Note that GN needs to choose the oversampling $l = 0.5r$ to guarantee its convergency. In the left plots, we show their relative errors to compare the approximation accuracy. In the right plots, the runtime of every calculation point is recorded to compare the computational speeds.

     In \ref{Fig.sub.sp} and \ref{Fig.sub.fp} of Fig.\ref{Fig.main.SP}, it is not surprised to find that GN-c performs the best in accuracy, and the GN-r\&c is slightly better than rSVD and GN. From the observations of \ref{Fig.sub.spt} and \ref{Fig.sub.fpt}, the computational speed of GN-c is slower than rSVD and GN, and GN is the fastest one.

     From the related statements in section 3 and the experimental results, GN-r\&c is not an efficient method because of its slow computational speed and much more storage requirement. GN-c outperforms the GN-r\&c and rSVD in accuracy. Generally speaking, the GN-c method is more advisable in finding a low-rank approximation for the nonsymmetric matrices with polynomially decaying spectrum.

     \subsubsection{Synthetic examples with exponentially decaying spectrum}
\begin{figure} [H]
    \subfigure[\textit{Synthetic matrix (ExpDecaySlow)}]{
        \includegraphics[width=8cm,height=6cm]{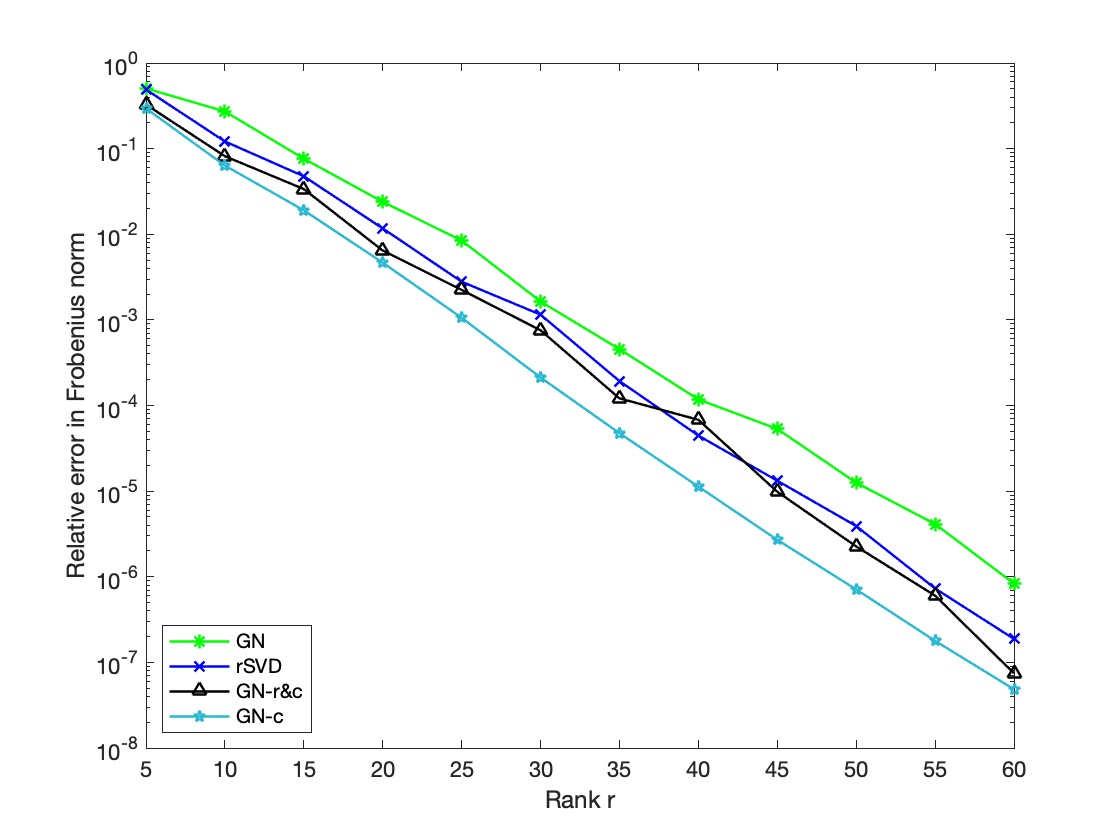} \label{Fig.sub.se}
    }
    \hspace{2mm}
    \subfigure[\textit{Synthetic matrix (ExpDecaySlow)}]{
        \includegraphics[width=8cm,height=6cm]{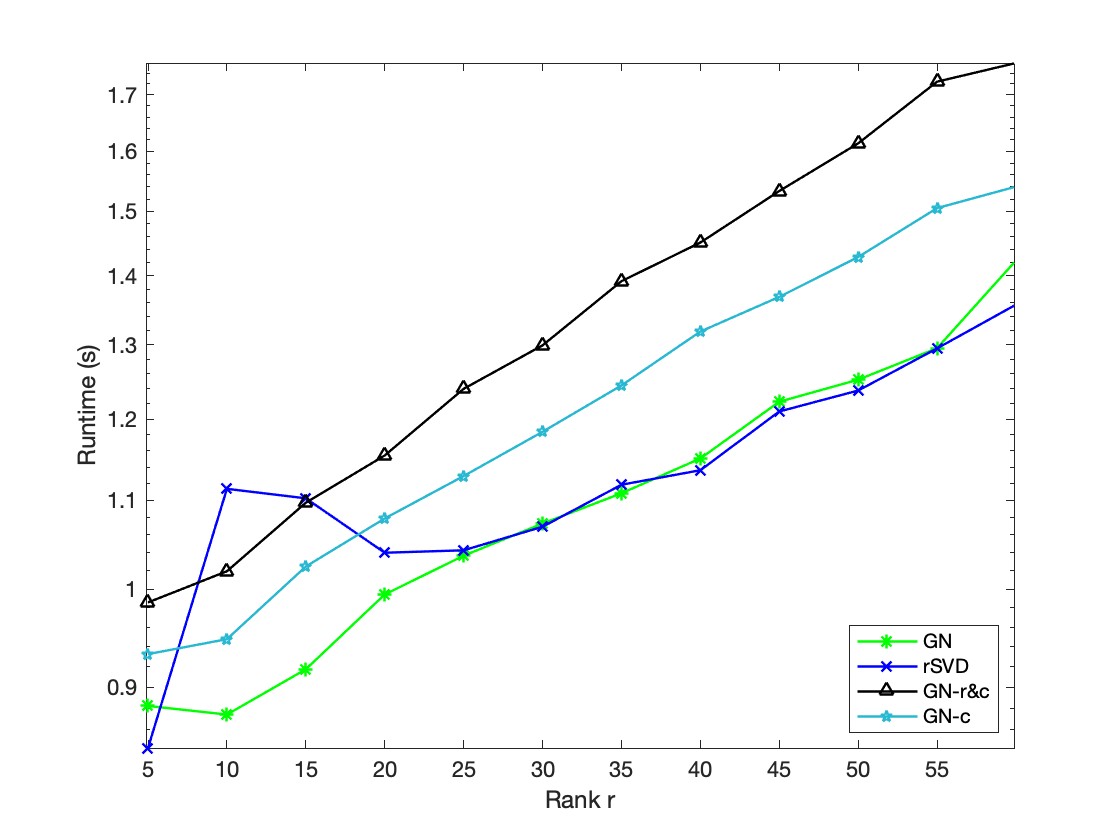} \label{Fig.sub.set}
    }

    \subfigure[\textit{Synthetic matrix (ExpDecayFast)}]{
        \includegraphics[width=8cm,height=6cm]{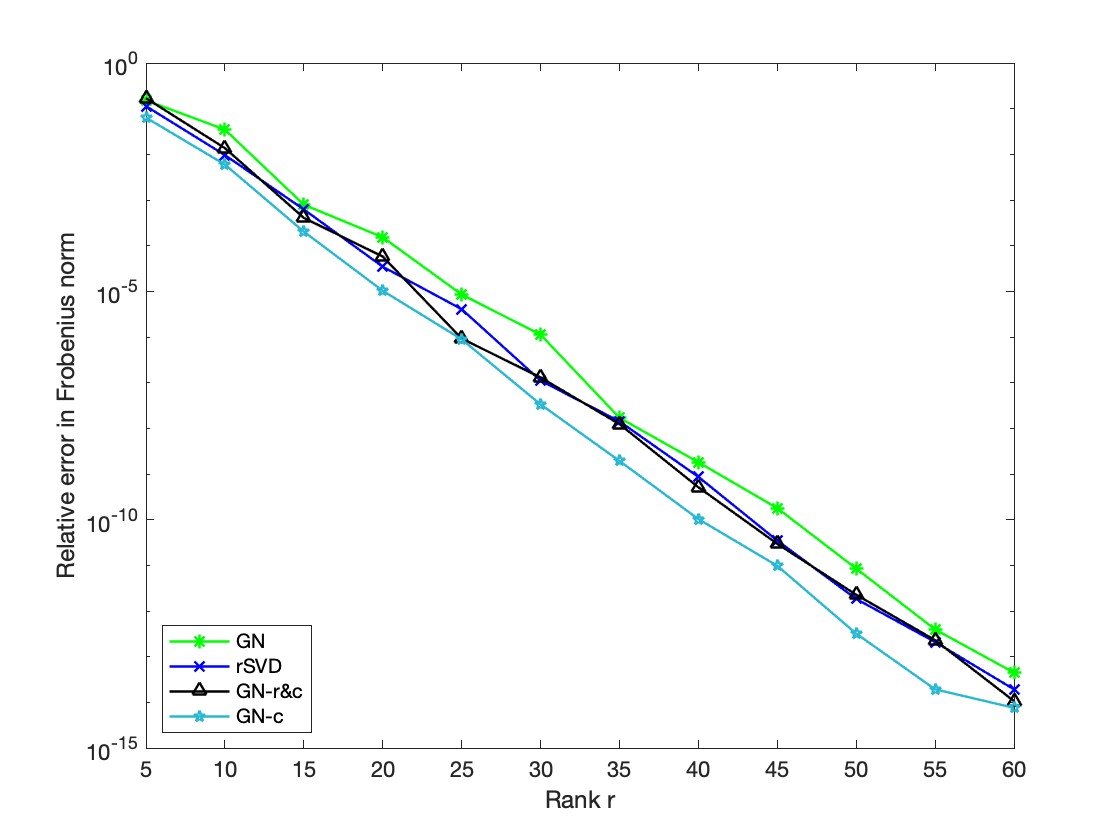} \label{Fig.sub.fe}
    }
    \hspace{2mm}
    \subfigure[\textit{Synthetic matrix (ExpDecayFast)}]{
        \includegraphics[width=8cm,height=6cm]{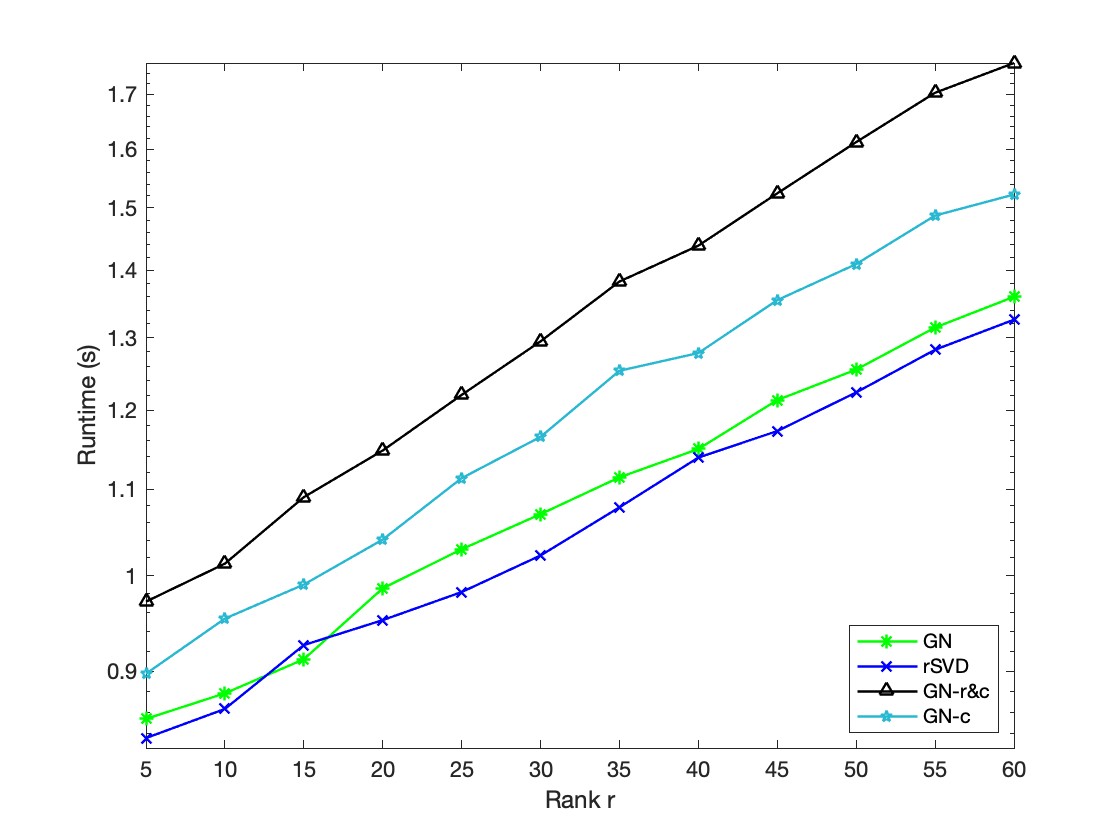} \label{Fig.sub.fet}
    }
    \caption{\textit{\textbf{Relative error and runtime comparisons.}  The relative error and runtime of approximating the matrix with slow exponential decaying speed are shown in (\romannumeral 1) and (\romannumeral 2). And the relative error and runtime of approximating the matrix with fast exponential decaying speed are presented in (\romannumeral 3) and (\romannumeral 4).}}\label{Fig.main.SE}
    \end {figure}

     Next, we explore to obtain a rank-$r$ approximation of $A\in\mathbb{R}^{15000\times 15000}$ with slow ($q = 0.125$) and fast ($q = 0.25$) exponential decaying speeds of singular spectrum in Fig.\ref{Fig.main.SE}. For GN, we take $l = 5$ in this numerical experiments. The relative errors of all the methods are shown in the left plots to compare the approximation accuracy, and in the right plots, the runtime of every calculation point is recorded to compare the speeds.

     As shown in \ref{Fig.sub.se} and \ref{Fig.sub.fe} of Fig.\ref{Fig.main.SE}, for the cases with slow and fast exponential decaying speeds, all the methods nearly obtain a good low-rank approximation. GN-r\&c and rSVD have subtle differences in relative errors. GN-c outperforms other methods in accuracy as expected.
     Nonetheless, the GN-c is slower in comparison with rSVD and GN in \ref{Fig.sub.set} and \ref{Fig.sub.fet}. However, considering the stability and accuracy of GN-c, it is still a promising method as shown above.

    Recalling from the numerical examples in Fig.\ref{Fig.main.2} of section 3, for synthetic dense matrices some preliminary conclusions can be drawn.
    \begin{itemize}
        \item[(1)]
        All the experiments shown here suggest that GN-c is more accurate and stable when approximating the dense nonsymmetric matrices regardless of how rapidly the singular values decay.
        \item[(2)]
        The computational speed of GN-c is slightly slower than rSVD and GN for the dense cases.
        \item[(3)]
        GN-r\&c is not an attractive method in practice. Though the GN-r\&c method tends to perform better than rSVD and GN in accuracy when the singular spectrum is polynomially decaying, it needs nearly double computation in comparison with GN-c.
    \end{itemize}

    Nextly we continue the numerical experiments to explore the low-rank approximation of the large-scale sparse dataset generated from various application scenarios.
    \subsection{Dataset examples}
   The large-scale nonsymmetric matrices in this subsection are selected at random from the matrix market mentioned earlier. Similarly, we still concentrate on the relative errors in the Frobenius norm and the runtimes to compare the performances of GN, rSVD, GN-r\&c, and GN-c.

\begin{figure} [H]
    \subfigure[\textit{fidap011 (Fluid dynamics)}]{
        \includegraphics[width=8cm,height=6cm]{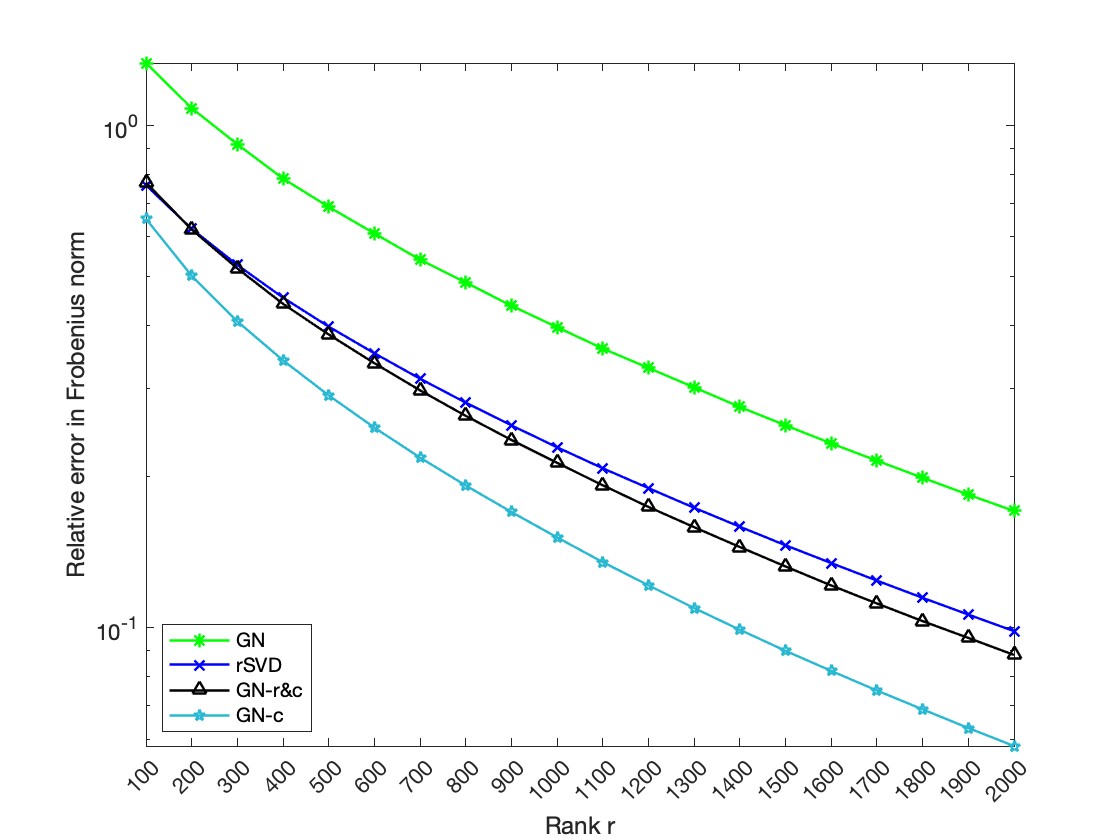} \label{Fig.sub.fidap}
    }
    \hspace{2mm}
    \subfigure[\textit{fidap011 (Fluid dynamics)}]{
        \includegraphics[width=8cm,height=6cm]{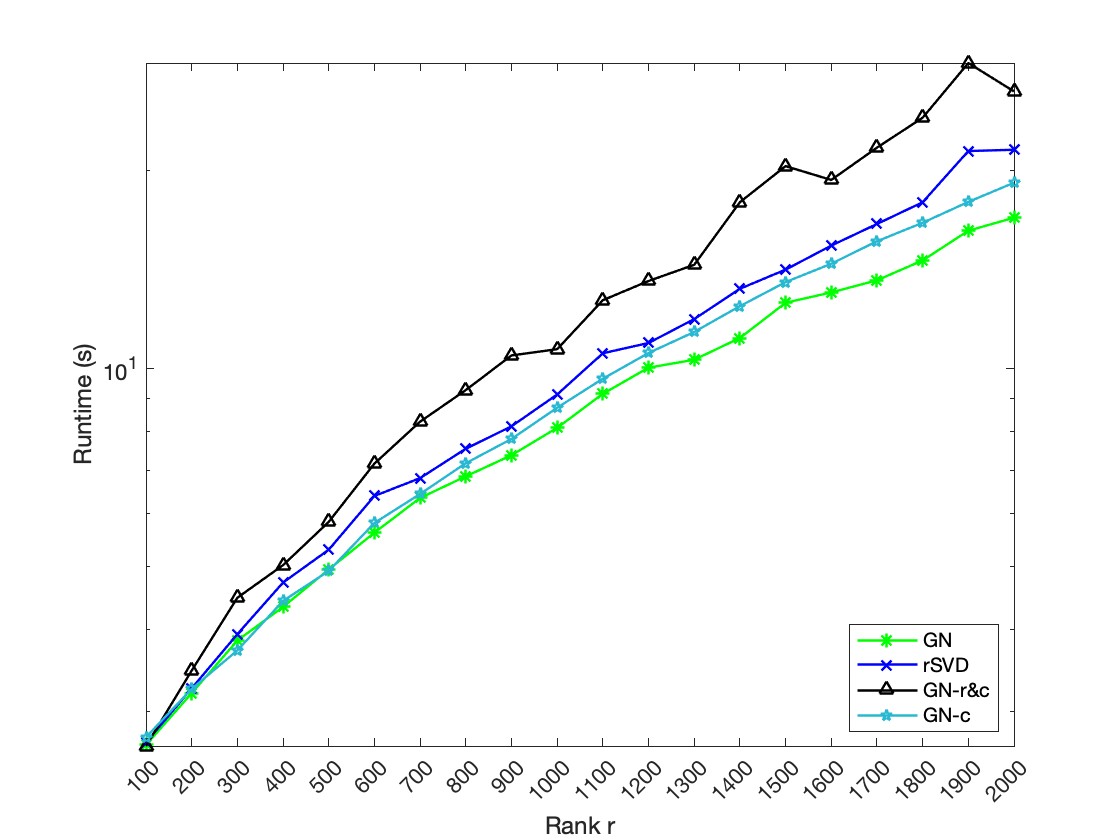} \label{Fig.sub.fidapt}
    }

    \subfigure[\textit{bcsstm25 (Structural problem)}]{
        \includegraphics[width=8cm,height=6cm]{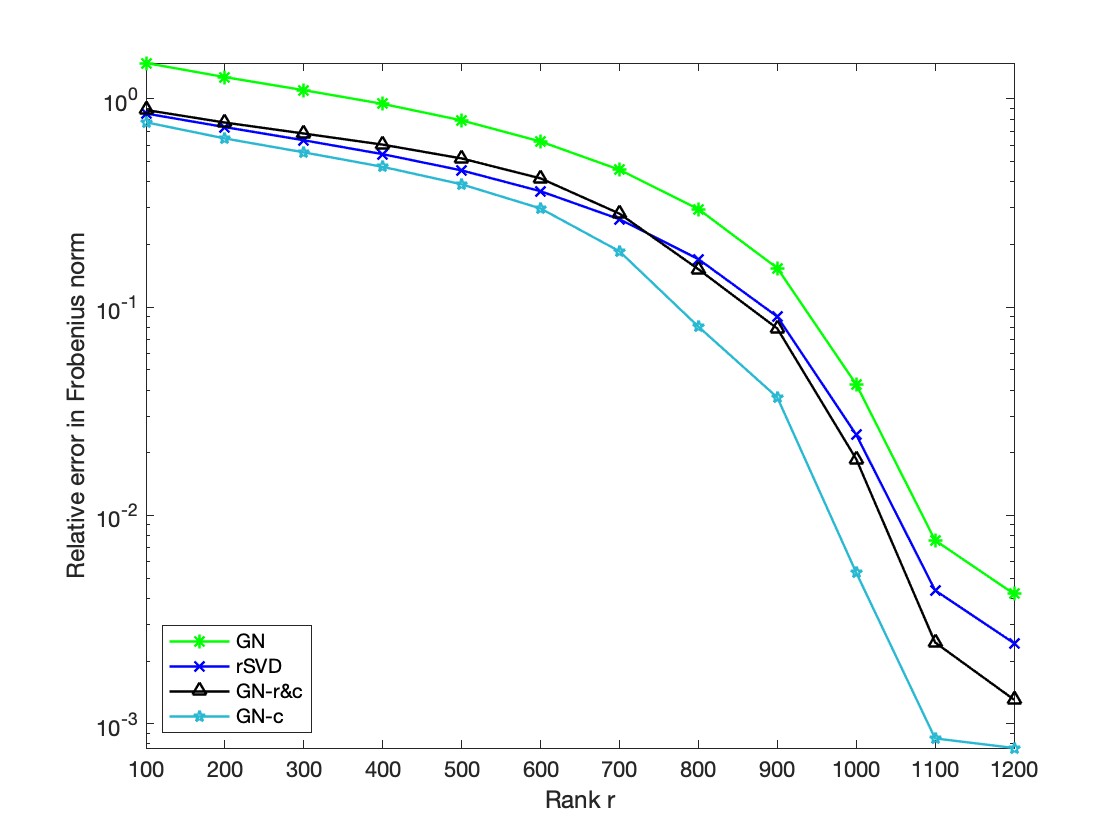} \label{Fig.sub.bcsstm}
    }
    \hspace{2mm}
    \subfigure[\textit{bcsstm25 (Structural problem)}]{
        \includegraphics[width=8cm,height=6cm]{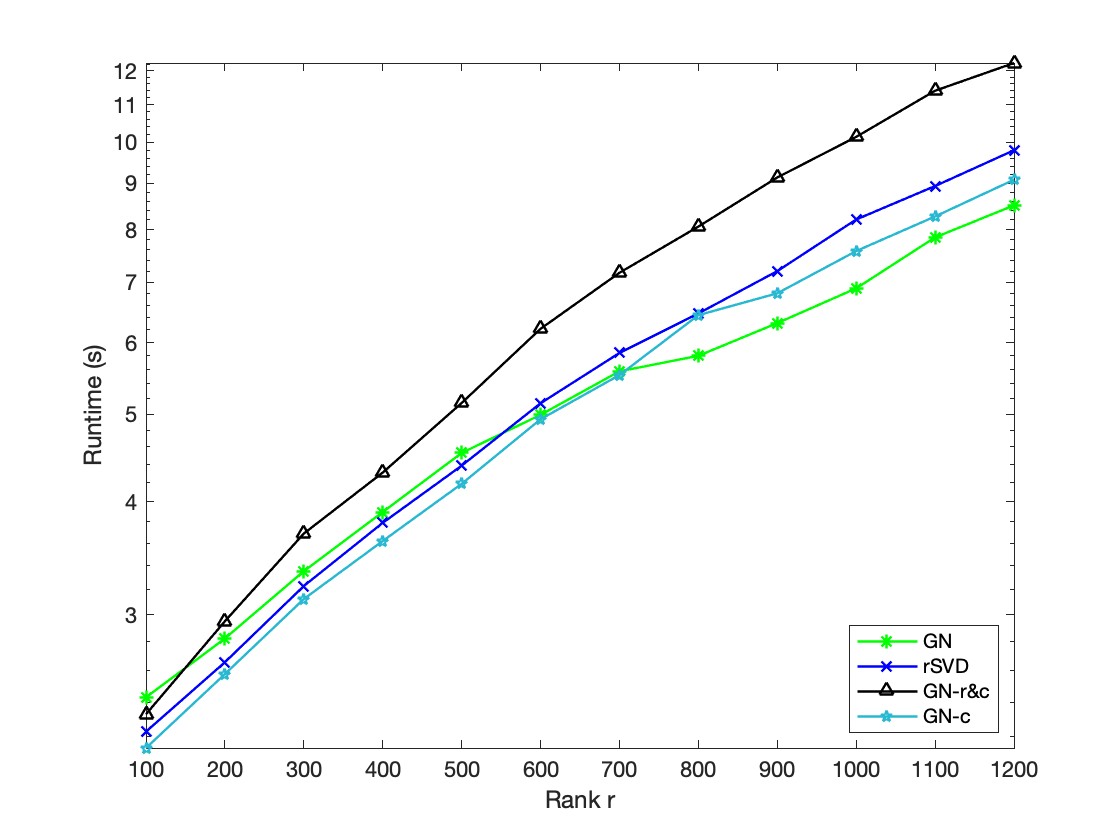} \label{Fig.sub.bcsstmt}
    }
    \caption{\textit{\textbf{Relative error and runtime comparisons.}  The relative error and runtime of approximating the matrix fidap011 are shown in (\romannumeral 1) and (\romannumeral 2). And the relative error and runtime of approximating the matrix bcsstm25 are presented in (\romannumeral 3) and (\romannumeral 4).}}\label{Fig.main.fidap}
    \end {figure}

   In Fig.\ref{Fig.main.fidap}, we test two sparse nonsymmetric matrices. The matrix fidap011 about fluid dynamics problem is $16614\times 16614$ and has $1091362$ nonzeros. And bcsstm25 is $15439\times 15439$ and owns 252241 nonzero entries, which is related to the structural problem.
    In \ref{Fig.sub.fidap} and \ref{Fig.sub.bcsstm}, GN-c performs much better than other methods in accuracy as expected. For the results in \ref{Fig.sub.fidapt} and \ref{Fig.sub.bcsstmt}, we can find that the computational speed of GN-c is faster than rSVD and GN-r\&c, and the performance of GN-c is still outstanding for the sparse cases.

\begin{figure} [H]
    \subfigure[\textit{jan99jac060 (Economic problem)}]{
        \includegraphics[width=8cm,height=6cm]{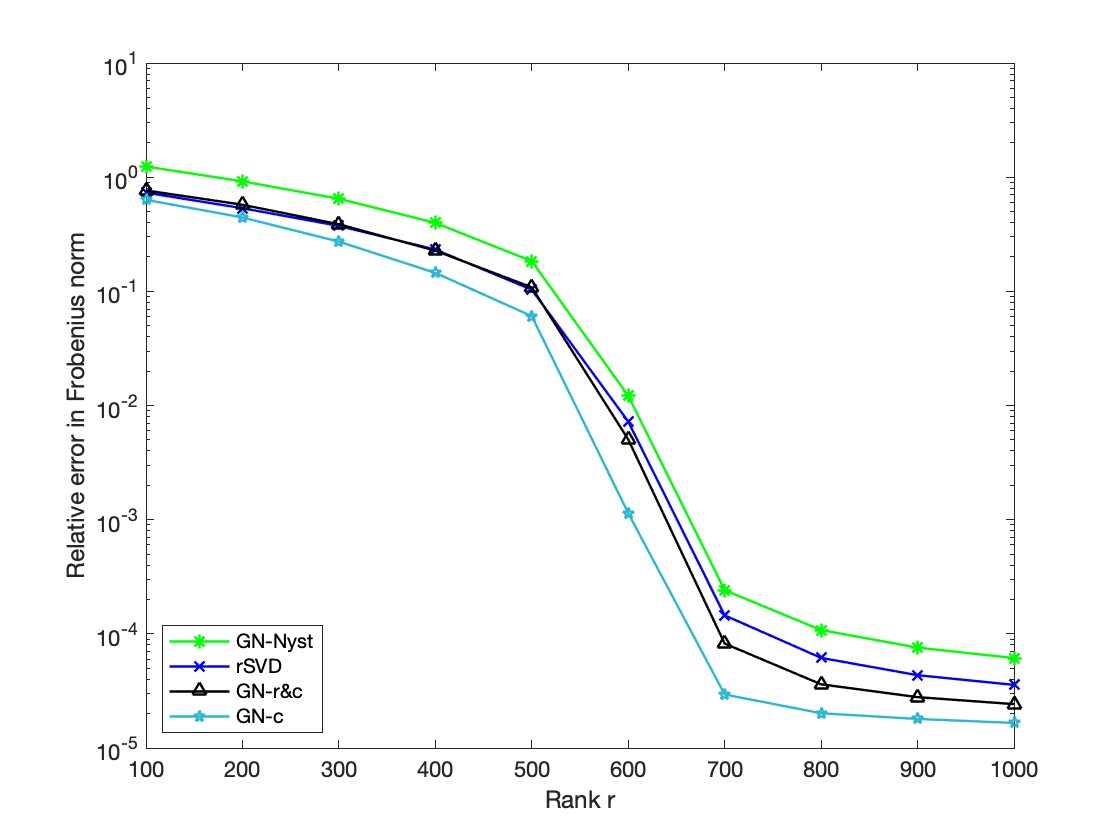} \label{Fig.sub.jan}
    }
    \hspace{2mm}
    \subfigure[\textit{jan99jac060 (Economic problem)}]{
        \includegraphics[width=8cm,height=6cm]{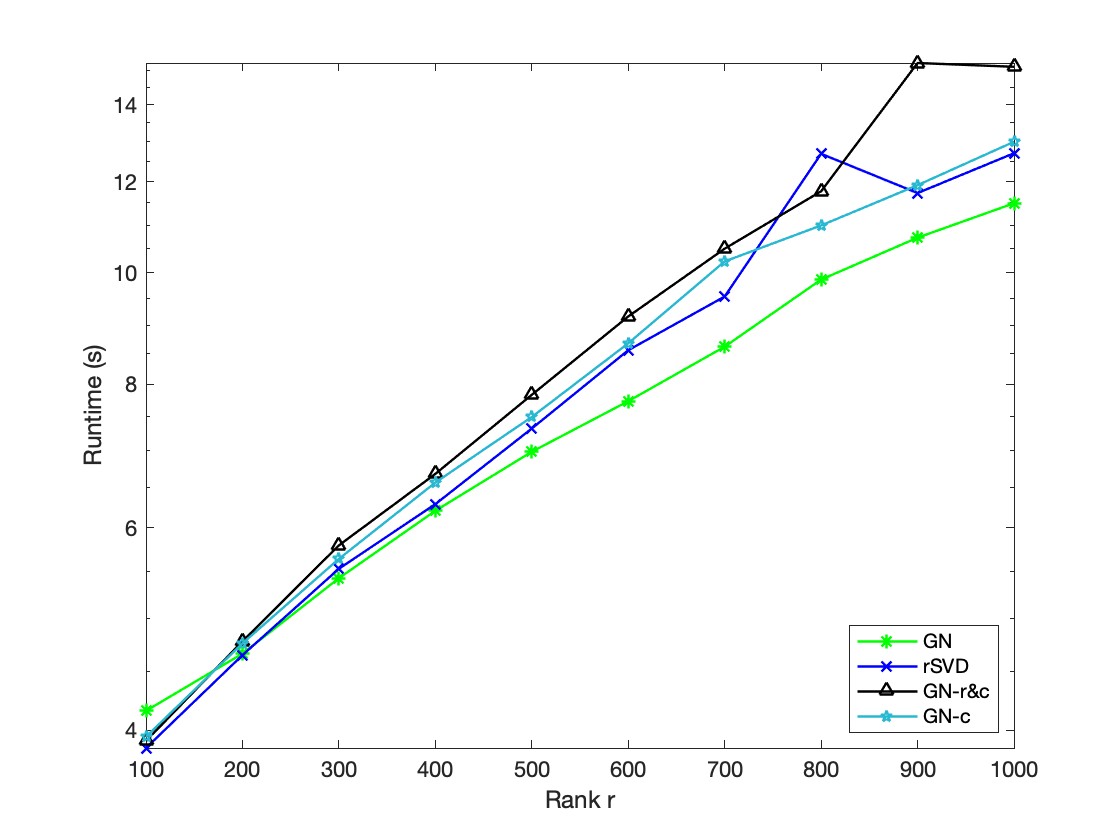} \label{Fig.sub.jant}
    }

     \subfigure[\textit{bayer10 (Chemical problem)}]{
        \includegraphics[width=8cm,height=6cm]{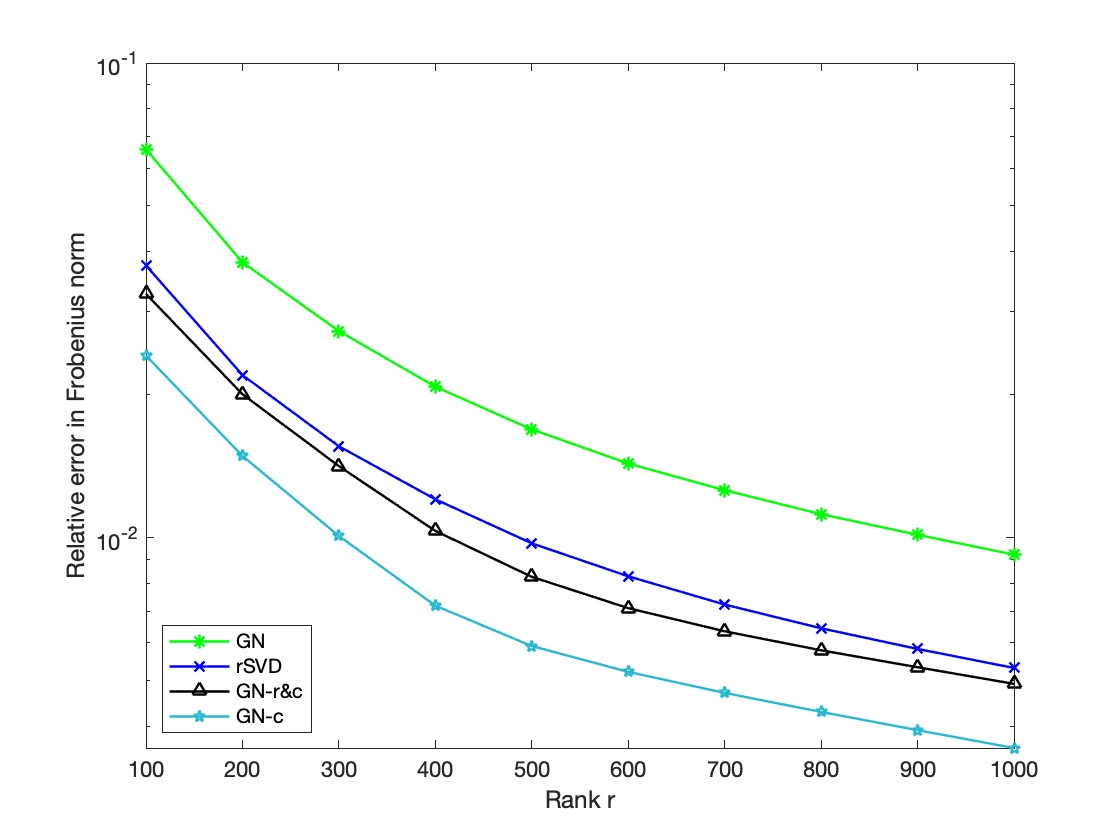} \label{Fig.sub.bayer}
     }
     \hspace{2mm}
     \subfigure[\textit{bayer10 (Chemical problem)}]{
        \includegraphics[width=8cm,height=6cm]{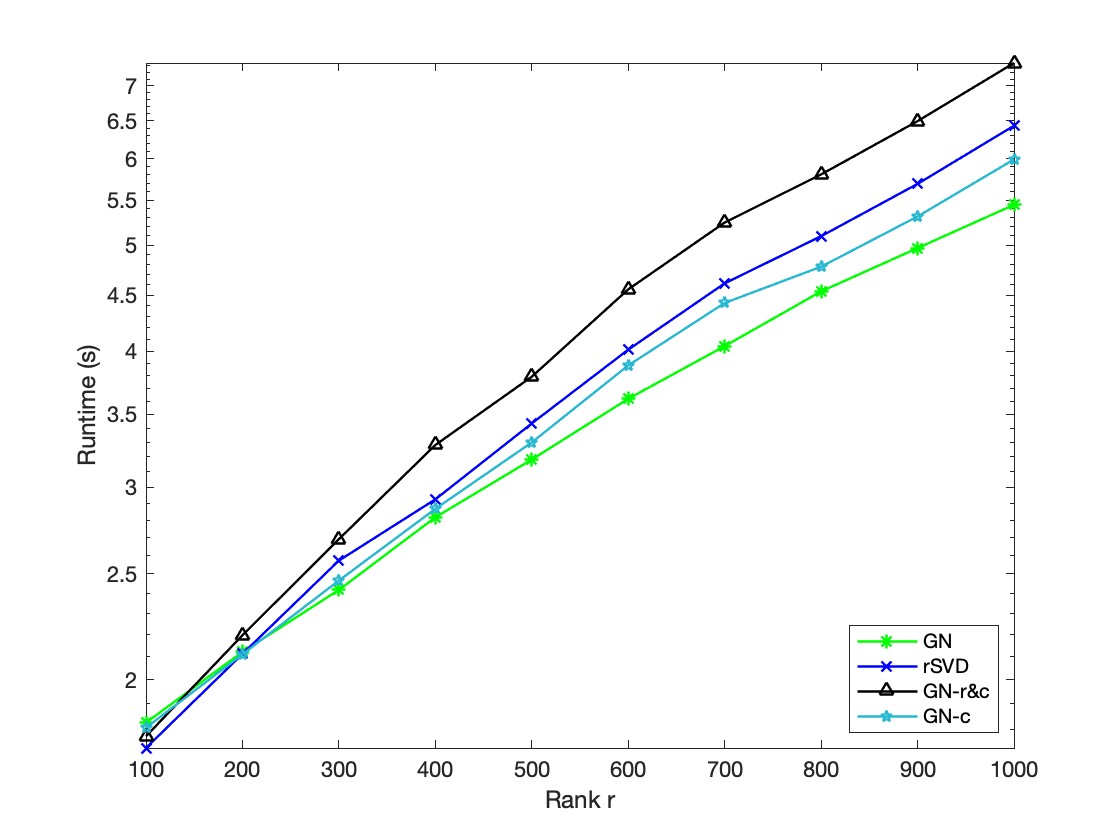} \label{Fig.sub.bayert}
     }
        \caption{\textit{\textbf{Relative error and runtime comparisons.}  The relative error and runtime of approximating the matrix jan99jac060 are shown in (\romannumeral 1) and (\romannumeral 2). And the relative error and runtime of approximating the matrix bayer10 are presented in (\romannumeral 3) and (\romannumeral 4).}}\label{Fig.main.jan}
    \end {figure}

     Furthermore, we choose another two sparse matrices jan99jac060 and bayer10.
     The matrix jan99jac060 associated with economic problem is $20614\times 20614$ and has 127182 nonzeros. The matrix bayer10 arises from chemical process simulation, which is $13436\times 13436$ and has 71594 nonzero entries.
     The \ref{Fig.sub.jan} and \ref{Fig.sub.bayer} in Fig.\ref{Fig.main.jan} show us that GN-c is notably more efficient in stability and accuracy.
     In \ref{Fig.sub.jant} and \ref{Fig.sub.bayert}, the speeds of GN-c and rSVD are almost the same and GN is the fastest as before.

     Obviously, all the experiments identify that GN-c outperforms GN and rSVD in accuracy and stability no matter for the synthetic dense matrices or the sparse matrices. As for the computational speed, GN-c is slightly faster than rSVD for some sparse cases, and slower than rSVD and GN when approximating the synthetic dense matrices. So generally GN-c is an efficient approach to find a low-rank approximation of nonsymmetric matrices in practice.
    \section{Conclusions}
    In this paper we have proposed a new algorithm GN-c for efficiently computing a low-rank approximation of nonsymmetric matrices, which is based on the framework of generalized Nystr\"{o}m form $A \approx AX(Y^{T}AX)^{\dagger}Y^{T}A$. The main techniques utilized by us are the rank-revealing factorizations and the randomized embeddings. Meanwhile, we have analyzed the error bound of GN-c in the Frobenius norm in section 4, which is identical with the results of the numerical experiments in section 5. Though not all methods can be expected to always perform well in practice, We have shown experimentally that the GN-c method outperforms GN and rSVD in accuracy, and verified that GN-c is an efficient method to obtain an accurate approximation for nonsymmetric matrices without sacrificing stability. Meanwhile in order to improve the performance of GN-c, we can explore to utilize other more efficient randomized rank-revealing QR factorizations. Additionally, the accuracy and computational speed of GN-c are theoretically associated with the decaying speed of singular spectrum. Many problems are worthwhile to explore ,especially when the singular spectrum of $A$ is flat.

\end{document}